
\magnification 1200
\input amstex
\documentstyle{amsppt}
\rightheadtext{Invariant tensors and normal forms}
\topmatter
\leftheadtext{Peter Ebenfelt}
\date{\number\year-\number\month-\number\day}\enddate


\define\Span{\text{\sl span }}

\define \im{\text{\rm Im }}
\define \re{\text{\rm Re }}

\define \bR{\Bbb R}

\define \bC{\Bbb C}


\def\a {\alpha}

\def\dim {\text {\rm dim}}

\def\Aut{\text{\rm Aut}}
\def\pb11{\overline{p_{11}}}
\def\pb02{\overline{p_{02}}}

\def\lpb{\left<}
\def\rpb{\right>}
\def\T{\Cal T}

\title  New invariant tensors in CR structures and a normal form for
real hypersurfaces at a generic Levi degeneracy\endtitle 
\author Peter Ebenfelt\footnote{Supported in part by a grant from the
Swedish Natural Science Research
Council.\newline}
\endauthor
\address Department of Mathematics, Royal Institute of Technology, 100
44 Stockholm, Sweden\endaddress
\email ebenfelt\@math.kth.se\endemail


\endtopmatter
\document

\heading 0. Introduction\endheading

In 1974, Chern and Moser [CM] solved the biholomorphic equivalence problem
for real-analytic hypersurfaces in $\bC^{n+1}$ at Levi nondegenerate
points. (The case $n=1$
was considered and solved by E. Cartan [C1--2].) They presented a
complete set of biholomorphic invariants for such a hypersurface
at a Levi nondegenerate point; by a complete
set of invariants, we mean a set of invariants such that given two
hypersurfaces $M,M'\subset\bC^{n+1}$ with distinguished points
$p_0\in M$, $p_0'\in M'$, there is a 
biholomorphic transformation $Z'=H(Z)$ near $p_0$ such that
$H(M)\subset M'$ and $H(p_0)=p_0'$ if and only if the set of
invariants for $M$ and $M'$ are equal. The Chern-Moser invariants can
in principle (there is an infinite number of invariants) be computed
from the Chern-Moser normal form, which is a normal form for a Levi
nondegenerate hypersurface $M$, defined in terms of the Levi form at
$p_0\in M$, such that the transformation to normal form
is unique modulo a finite dimensional normalization. 

In the present paper, we introduce a new sequence of invariant
tensors, $\psi_2,\psi_3\ldots$, 
for generic submanifolds of $\bC^N$ (Theorem 2.9), which can be viewed
as higher order Levi 
forms. (Although the tensors are only introduced here in the context of
generic submanifolds of $\bC^N$, it is clear that the definitions
work equally well in general CR structures.) The second order tensor $\psi_2$
coincides with the Levi map and the higher order tensors are related,
as explained in \S3 below, 
to the data of {\it finite nondegeneracy}, a notion which has recently
proved very useful in the study of real submanifolds in $\bC^N$ (see
e.g.\ [BER1--4], [BHR], [E1--2]). The third order tensor is also
related to the {\it cubic form} as introduced by Webster [W] (see
Remark 4.17). 

As one of our main results (Theorem 1.1.28), we describe, 
using the second and third order tensors, a formal normal
form (in the sense of Chern--Moser as described above) for a real
smooth (meaning $C^\infty$) hypersurface  
$M\subset\bC^{n+1}$ at a {\it generic Levi degeneracy} $p_0\in M$,
i.e. a point $p_0$ at which the Levi determinant  
vanishes but its differential does not and the set of Levi
degenerate points of $M$ (which is then a smooth codimension one submanifold of
$M$ at $p_0$) is transverse to the Levi null space (which is then one
dimensional) at that point. (We refer the reader to [W] for further
discussion of generic Levi degeneracies; for instance, a normal form for generic
Levi degeneracies in $\bC^2$ 
under formal holomorphic {\it contact} transformations is given in [W].) In view of a 
convergence result due to the author, 
Baouendi, and Rothschild [BER3], the formal normal form  in Theorem 1.1.28
provides a complete set of {\it biholomorphic} invariants if the
hypersurface is also real-analytic (Corollary
1.1.30). 

We then proceed to study the special case where the Levi
form, at the generic Levi degeneracy, is semidefinite. In this
situation, the normal form can be expressed in a
particularly simple 
and explicit form (Theorem 1.2.5) by applying a result of E. Cartan. (The associated
partial (third order) normal form is given with numerical invariants;
in fact, an explicit partial normal form which is valid in a slightly
more general setting is given in Theorem 1.2.10.) The corresponding explicit
character of the normalization of the transformation to normal form
makes it 
possible to compute bounds for the stability group of a real
hypersurface at a generic semidefinite Levi degeneracy (Corollary
1.2.7). In the case $n=2$, i.e.\ for
hypersurfaces in $\bC^3$, the results on normal forms in this paper
are contained in the results of [E1]. However, the invariant tensors
introduced here shed additional light on some of the results from, and
answers a question posed in, that paper. 

The paper is organized as follows. In the first section, \S1.1, we
present the normal form for a generic Levi degeneracy. In \S1.2, we
give the more explicit normal form in the special 
case where the Levi form at $p_0$ is also semidefinite. 
We then turn to the more general 
situation of generic submanifolds in $\bC^N$ and introduce the CR
invariant tensors. \S3 is devoted to
explaining the relation between the notion of finite nondegeneracy and
the tensors of \S2. In \S4, we return to the case of
hypersurfaces and show, as a preparation for the normal form, that the
second and third order tensors form a complete set of third order
invariants for a real hypersurface by relating these tensors to the
defining equation of $M$. Then, we calculate, in \S5, explicit numerical
invariants associated with the third order tensor of a real
hypersurface at a point where the Levi form has rank $n-1$ and is
semidefinite. \S6--8 are devoted to the proofs of the results that give
the normal form.

\subhead Acknowledgments\endsubhead The author would like to thank
M. S. Baouendi, L. P. Rothschild, and N. Wallach for valuable comments
on a preliminary version of this paper.

\heading 1. Normal forms for real hypersurfaces at generic Levi
degeneracies\endheading

\subhead 1.1. The general case\endsubhead In this section, we shall
present a normal form for a generic Levi 
degeneracy; the reader should recall the definition of generic Levi
degeneracy from \S0. In order to describe the
normal form, we need first a partial (third order) normal form. We
begin with some notation. 

We use the notation $\Cal M(\bC^m)$ for the space of $m\times m$
matrices with complex matrix elements and $GL(\bC^m)$ for
the group of invertible ones. We also write $\Cal S(\bC^m)$ for
the symmetric matrices in $\Cal M(\bC^m)$, i.e. those for which
$A=A^\tau$. (Here, $A^\tau$ denotes the transpose of $A$.) For
nonnegative integers $r$ and $s$ 
such that $r+s=m$, we denote by $\hat U(r,s,\bC)$ the subgroup of
$GL(\bC^m)$ consisting of those matrices $U$ for which
$$
U^* I_{r,s}U=\pm I_{r,s},\tag1.1.1
$$ 
where $I_{r,s}\in GL(\bC^m)$ is the diagonal matrix whose $r$ first 
diagonal elements are $+1$ and $s$ last ones are $-1$, and  $U^*$
denotes the Hermitian adjoint of $U$ (i.e. $U^*=\bar U^\tau$). This
group decomposes naturally 
as 
$$
\hat U(r,s,\bC)=U^+(r,s,\bC)\cup U^-(r,s,\bC), 
$$
where $U^+(r,s,\bC)$ and 
$U^-(r,s,\bC)$ denote the set of matrices for which \thetag{1.1.1}
holds with the $+$ and $-$ sign, respectively. Observe that
$U^+(r,s,\bC)$ is a subgroup whereas $U^-(r,s,\bC)$ is not. Also, note
that $U^-(m,0,\bC)$ is
empty, and $U(m,0,\bC)=U^+(m,0,\bC)$ is the usual group $U(\bC^m)$ of unitary
matrices. 

Consider the action of the group $\bR_+\times \hat U(r,s,\bC)$ on
$\Cal S(\bC^m)$ given by
$$
\bR_+\times \hat U(r,s,\bC)\ni (\sigma,U)\to (\sqrt{\sigma}U)^\tau A
(\sqrt{\sigma}U)\in \Cal S(\bC^m),\tag1.1.2
$$
for $A\in \Cal S(\bC^m)$. We denote, for given $A\in\Cal S(\bC^m)$, by
$C_{r,s}(A)\subset\Cal S(\bC^m)$ its orbit or conjugacy class under
the group action \thetag{1.1.2}, i.e. 
$$
C_{r,s}(A)=\left\{B\in\Cal S(\bC^m)\:B=(\sqrt{\sigma}U)^\tau A(\sqrt{\sigma}U),\ U\in
\hat U(r,s,\bC),\ \sigma>0\right\}.
$$
These conjugacy classes form a disjoint partition of $\Cal S(\bC^m)$. We
have the following result, which is 
the first step in describing the normal form and whose proof will be
given in \S6.

\proclaim{Proposition 1.1.3} Let $M\subset\bC^{n+1}$ be a real smooth
hypersurface. Assume that $M$ has a generic Levi degeneracy at $p_0\in
M$. Denote by $r$, with $(n-1)/2\leq r\leq n-1$, the number of eigenvalues
of the Levi form at $p_0$ which have the same sign. (Also, write
$s=n-1-r$.) Then, there exists
a unique conjugacy class 
$C_{r,s}\subset\Cal S(\bC^{n-1})$ and for each $R\in C_{r,s}$
there are local holomorphic
coordinates $Z=(z,w)\in\bC^n\times \bC$,
$z=(z',z^n)\in\bC^{n-2}\times\bC$, near $p_0$, vanishing at 
$p_0$, such that the defining equation of $M$ is of 
the following form,
$$
\im
w=\sum_{j=1}^{r}|z^j|^2-\sum_{k=r+1}^{n-1} |z^k|^2+
2\re\!\!\left(\bar
z^{\bar n}\left((z')^\tau Rz'+(z^n)^2\right)\right)+F(z,\bar
z,\re w),
\tag1.1.4
$$
where
$F(z,\bar z,\re w)$ denotes a smooth, real valued function
which is $O(4)$ in the weighted coordinate system where $z$ has weight
one and $w$ weight two. 
\endproclaim

Let us briefly explain our usage of the notation $O(\nu)$, for
nonnegative integers $\nu$, in Proposition 1.1.3. We assign the weight one
to the variables $z=(z',z^n)=(z^1,\ldots,z^{n-1},z^n)$, the weight  
two to $w$, and say that a polynomial $p_\nu(z,w)$ is weighted
homogeneous of degree $\nu$ if, for all $t>0$,
$$
p_\nu(tz,t^2w)=t^\nu p_\nu(z,w).\tag1.1.5
$$
We shall write $O(\nu)$ for a formal series involving only terms of
weighted degree greater than or equal to 
$\nu$. We say that a smooth function defined near $0$ is $O(\nu)$ (at
$0$) if its Taylor series at $0$ is $O(\nu)$. Similarly, we speak of
weighted homogeneity 
of degree $\nu$ and 
$O(\nu)$ for polynomials, power series, and functions in $(z,\bar z,\re w)$, where
$\bar z$ is assigned the weight one and $\re w$ the weight two.

We shall now present a complete, formal, normal form for a generic Levi
degeneracy. Before stating the 
theorem, we need to define the space of normal forms and the
normalization for the transformation to normal form.
 
By Proposition 1.1.3, we may
assume that $M$ is defined near $p_0=(0,0)$ by \thetag{1.1.4} for
given, and fixed for the remainder of this section, integer $r$ and matrix $R\in\Cal
S(\bC^{n-1})$. Since we
shall present a formal normal form, we consider the defining equation
\thetag{1.1.4} as a formal power series. It is well known (cf. [BJT] and
[BR]; cf. also the forthcoming book [BER4]) that, after an additional
formal change of coordinates at $(0,0)$ if necessary, we may also
assume that $F(z,0,s)\equiv F(0,\bar z,s)\equiv 0$;  we shall say that
the (formal) coordinates $(z,w)$ are {\it regular} for $M$ at
$p_0=(0,0)$ if the (formal) defining equation for $M$ at that point is
of the form $\im w=\phi(z,\bar z,\re w)$ with
$\phi(z,0,s)\equiv\phi(0,\bar z,s)\equiv 0$. We subject the (formal)
hypersurface $M$ to a formal 
invertible transformation
$$
z=\tilde f(\tilde z,\tilde w)\quad,\quad w=\tilde g(\tilde z,\tilde w),\tag1.1.6
$$
where 
$$
\tilde f=(\tilde f',\tilde f^n)=(\tilde f^1,\ldots ,\tilde
f^{n-1}, \tilde f^n),\tag1.1.7
$$ 
such that the form
\thetag{1.1.4} is preserved. We shall also require that the coordinates
$(\tilde z,\tilde w)$ are regular for $M$, i.e. the remainder $\tilde
F(\tilde z,\bar{ \tilde z},\tilde s)$ corresponding to the defining
equation relative the coordinates $(\tilde z,\tilde w)$ satisfies
$\tilde F(\tilde z,0,\tilde s)\equiv\tilde F(0,\bar{\tilde z},\tilde
s)\equiv 0$. 

Given a matrix $A\in\Cal S(\bC^m)$, we denote by
$O_A(\bC^m)$ the subgroup of $GL(\bC^m)$ consisting of those matrices
that preserve the bilinear form associated with $A$, i.e. $B\in
O_A(\bC^m)$ if 
$$
B^\tau AB=A.\tag1.1.8
$$

We have the following proposition whose proof will be given in \S7.
\proclaim{Proposition 1.1.9} A transformation \thetag{1.1.6} preserving
regular coordinates also preserves the
form \thetag{1.1.4}, for a given integer $r$ and $R\in\Cal S(\bC^{n-1})$,
if and only if
$$
\aligned
\tilde f'(z,w)& =Az'+wB+\frac{2i}{c}(B^* Az')Az'+O(3)\\
\tilde f^n(z,w)& =c^{1/3}z^n+O(2)\\
\tilde g(z,w)& =cw+2i(B^* Az')w+O(4),\endaligned\tag1.1.10
$$
where $c\in\bR\setminus\{0\}$, $B\in\bC^{n-1}$ (considered as an $(n-1)\times 1$
matrix), and where $A$ is such that $A/|c|^{1/3}\in  O_R(\bC^{n-1})$ and 
$A^*I_{r,s}A=cI_{r,s}$ (in particular then, $A/|c|^{1/2}\in \hat
U(r,s,\bC)$); As above, we use $s=n-r-1$.
\endproclaim
 
We shall consider formal mappings \thetag{1.1.6} of the following form
$$
(\tilde f(z,w),\tilde g(z,w))=(T\circ P)(z,w).\tag 1.1.12
$$
Here, $P(z,w)$ is a polynomial mapping, $P=(P',P^n,P^{n+1})$,
$$
P(z,w)=(p'(z,w)+q'(z,w),p^n(z,w)+q^n(z,w),p^{n+1}(z,w)),\tag1.1.13
$$
where $p'=(p^1,\ldots, p^{n-1})$, $p^n$, $p^{n+1}$ are polynomials of 
the form
$$
\aligned
p'(z,w)& =Az'+wB+\frac{2i}{c}(B^* Az')Az'\\
p^n(z,w)& =c^{1/3}z^n\\
p^{n+1}(z,w)& =cw+2i(B^* Az')w,\endaligned\tag1.1.14
$$
where $c$, $B$, $A$  are as in Proposition 1.1.9. 
The polynomials $q'=(q^{1},\ldots, q^{n-1})$ and $q^n$ are
weighted homogeneous of the forms
$$
\aligned
q^\beta(z,w) &=\sum_{|J|=3} a^\beta_J
z^J+\left(\sum_{\alpha<\beta}b^\beta_\alpha (Az')^\alpha+c^\beta(Az')^\beta\right)w,\\
q^{n}(z,w) &=\sum_{|I|=2}d_I z^I,\endaligned\tag1.1.15
$$
where $\beta=1,\ldots, n-1$, $a^\beta_J, b^\beta_J,d_I\in\bC$,  and
$c^\beta\in\bR$. We use here multi-index notation so that
e.g. $J=(J_1,\ldots, J_n)$, $|J|=\sum_k J_k$, and
$z^J=(z^1)^{J_1}\ldots (z^n)^{J_n}$. The notation $(Az')^\beta$ stands
for the $\beta$:th component of the vector $Az'$. $T(z,w)$
in \thetag{1.1.12} is a formal mapping of the form
$$
T(z,w)=(z+f(z,w),w+g(z,w)),\tag1.1.16
$$
where $f=(f',f^n)=(f^1,\ldots, f^{n-1},f^n)$, and $g$ are formal power series in $(z,w)$ such
that $f'$ is $O(3)$, $f^n$ is $O(2)$, and $g$ is $O(4)$. We shall also
require that the formal
series $f'$, $f^n$ are such that the constant terms in the following
formal series vanish
$$
\frac{\partial^2 f^n}{\partial z^I}\,,\,  \frac{\partial^3
f^\beta}{\partial z^J} \,,\,
\re
\frac{\partial^2
f^\beta}{\partial z^\beta\partial w}\,,\,\frac{\partial^2
f^\beta}{\partial z^\alpha\partial w},\tag1.1.17
$$
where $I$ and $J$ range over all the multi-indices with $|I|=2$ and
$|J|=3$, respectively, the index $\beta$ runs over $1,\ldots, n-1$, and
$\alpha$ runs over $1,\ldots, \beta-1$. 
It is straightforward, and left to the reader, to verify (using
Proposition 1.1.9) that any
formal mapping \thetag{1.1.6} that preserves the form \thetag{1.1.4} of $M$
can be factored uniquely according to \thetag{1.1.12} with $T$ and $P$ as
above. We shall say that a choice of $P$, as described above, is a choice of
{\it normalization} for the transformations that preserve the form
\thetag{1.1.4} and that a formal mapping preserving the form
has this normalization if it is factored according to \thetag{1.1.12}
with this $P$.

Now, let $F(z,\bar z,s)$ be a
formal series in $(z,\bar z,s)$. In what follows,
we shall decompose the formal series $F(z,\bar z,s)$ as follows,
$$
F(z,\bar z, s)=\sum_{k,l}F_{kl}(z,\bar z,s),\tag1.1.18
$$
where $F_{kl}(z,\bar z,s)$ is of type $(k,l)$ i.e. for each $t_1,t_2>0$
$$
F_{kl}(t_1z,t_2\bar z,s)=t_1^kt_2^lF_{kl}(z,\bar z,s).\tag1.1.19
$$
We shall consider only those $F(z,\bar z,s)$ which are $O(4)$ and
which are ``real'' in the sense that
$$
F_{kl}(z,\bar z,s)=\overline{F_{lk}(z,\bar z,s)}.\tag1.1.20
$$
We shall denote by $\Cal F$ the space of all such formal power
series, and by $\Cal F_{kl}$ the space consisting of those which have
type $(k,l)$. In what follows, $F_{kl},H_{kl}$, and $N_{kl}$ denote formal
power series in $\Cal F_{kl}$.

In order to describe the space of normal forms, $\Cal N\subset\Cal F$,
we need a little more notation. Recall that the integer $r$ and matrix
$R\in\Cal S(\bC^{n-1})$ from Proposition 1.1.3 are fixed throughout this
section. For $u=(u^1,\ldots,u^{n-1})$ and
$v=(v^1,\ldots, v^{n-1})$, we use the notation
$\left<\cdot,\cdot\right>$ for the bilinear form
$$
\left<u,v\right>=\sum_{j=1}^{r}u^jv^j-\sum_{k=r+1}^{r}u^kv^k.\tag1.1.21
$$
We denote by $p_{R}(z)$ the quadratic polynomial
$$
p_{R}(z)=(z')^\tau Rz'+(z^n)^2.\tag1.1.22
$$
We use the notation
$\nabla=(\nabla',\nabla_n)=(\nabla_1,\ldots,\nabla_{n-1},\nabla_n)$
for the holomorphic gradient
$$
\nabla=\left(\frac{\partial}{\partial
z^1},\ldots,\frac{\partial}{\partial z^n}\right), 
$$
and similarly for the anti-holomorphic gradient $\bar \nabla$. We
shall need the linear operator $S_R$ defined on formal series
$u=u(z,\bar z,s)$ as follows
$$
S_Ru=-\left<\nabla',\bar\nabla'\right>(p_{R}u).\tag1.1.23
$$
Observe that $S_R$ maps $\Cal F_{k-1,l+1}$ into $\Cal F_{kl}$. Let us
remark that the operator $\left<\nabla',\bar\nabla'\right>$ is
essentially the same as the contraction operator tr corresponding to
the bilinear form $\left<\cdot,\cdot\right>$ as defined in [CM]; they
correspond to different normalizations for the monomials. 

We define the space of normal forms
$\Cal N\subset\Cal F$ for $M$ of the form \thetag{1.1.4} with $r$ and
$R$ as in Proposition 1.1.3 as follows. First, a formal series 
$N(z,\bar z,s)$ in $\Cal N$ is in {\it regular form} which can be expressed by 
$$
N(z,\bar z,s)=\sum_{\min(k,l)\geq 1}N_{kl}(z,\bar z,s);\tag1.1.24
$$
thus, $N$ has no components of type $(k,l)$ with $k=0$ or $l=0$. 
Moreover, the nonzero terms $N_{kl}$ satisfy the following
conditions:
$$
\aligned
N_{22}\in\Cal N_{22},&\quad N_{32}\in \Cal N_{32}\\
N_{42}\in\Cal N_{42},&\quad N_{33}\in \Cal N_{33}\\
N_{k1}\in\Cal N_{k1},&\quad k=1,2,3\ldots,\endaligned
\tag1.1.25
$$
where
$$
\aligned
N_{11} =& \big\{F_{11}\:F_{11}\in\ker \left<\nabla',\bar\nabla'\right>\big\}\\
N_{22} =& \big\{F_{22}\:F_{22}=\left<z',\bar z'\right>z^n\bar z^{\bar
n} H_{00}+H_{22},\
H_{22}\in \ker \left<\nabla',\bar\nabla'\right>\big\}\\
N_{33} =& \big\{F_{33}\:F_{33}=\left<z',\bar z'\right>^2(z^n
H_{01}+\overline{z^n H_{01}})+H_{33},\
H_{33}\in \ker
\left<\nabla',\bar\nabla'\right>^2\big\}\\
N_{21} =& \big\{F_{21}\:F_{21}=\bar z_{\bar n}H_{20}\big\}\\
N_{31} =&\big\{F_{31}\:F_{31}\in\ker S_R\big\}\\
N_{32} =&\big\{F_{32}\:F_{32}=\left<z',\bar
z'\right>^2z^nH_{00}+\left<z',\bar z'\right>H_{21}+H_{32},\\&
\hskip 1.15truecm
H_{32}\in\ker\left<\nabla',\bar\nabla'\right>,\
H_{21}\in\ker p_R(\nabla),\
\left<\nabla',\bar\nabla'\right>H_{21}\in\im S_R\big\}
\\
N_{42} =&\big\{F_{42}\:F_{42}=\left<z',\bar z'\right>\bar z^{\bar
n}H_{30}+H_{42},\ H_{42}\in\ker \left<\nabla',\bar\nabla'\right>,\\&
\hskip 1.15truecm
H_{30}\in\ker\nabla_n\big\},\endaligned\tag1.1.26
$$
and finally, for $k\geq4$,
$$
N_{k1} =\big\{F_{k1}\:F_{k1}=\bar z^{\bar n}H_{k0},\
H_{k0}\in\ker\nabla_n\big\}\tag1.1.27
$$
Observe that, for a series $H_{k0}$ of type $(k,0)$,  the
condition $H_{k0}\in\ker\nabla_n$ is equivalent to the condition that
$H_{k0}$ is independent of $z^n$, i.e. $H_{k0}=H_{k0}(z',s)$. 

We are now in a position to state the theorem on normal forms for a
generic Levi degeneracy.

\proclaim{Theorem 1.1.28} Let $M$ be a smooth hypersurface in
$\bC^{n+1}$ given
near $0\in M$ by \thetag{1.1.4}, where $r$ and $R$ are as in
Proposition $1.1.3$. Then, given any choice of
normalization (i.e. a choice of $P$ as described above), there is a
unique formal 
transformation \thetag{1.1.6} with this normalization that transforms
the defining equation \thetag{1.1.4} of $M$ at $0$ to
$$
\multline
\im
w=\sum_{j=1}^{r}|z^j|^2-\sum_{k=r+1}^{n-1} |z^k|^2+
2\re\!\!\left(\bar
z^{\bar n}\left((z')^\tau Rz'+(z^n)^2\right)\right)\\+N(z,\bar
z,\re w),\endmultline
\tag1.1.29
$$
where $N(z,\bar z,s)\in\Cal N$. \endproclaim

The proof of Theorem 1.1.28 will be given in \S8. We conclude this
section by applying Theorem 1.1.28 to the biholomorphic equivalence
problem. Suppose that 
$(M,p_0)$ and $(M',p_0')$ 
are two germs of real-analytic hypersurfaces in $\bC^{n+1}$ which have
generic Levi degeneracies at $p_0$
and $p_0'$, respectively. 
Thus, $M$ and $M'$ are, in particular, finitely nondegenerate (see
\S3) at their
distinguished points $p_0$ and $p_0'$. In view of [BER3, Theorem 2.6], any formal
equivalence between $(M,p_0)$ and $(M',p_0')$ is then in fact
biholomorphic. Hence, an immediate consequence of Theorem 1.1.28, as in
[E1], is the following.

\proclaim{Corollary 1.1.30} Let $M$ and $M'$ be 
real-analytic hypersurfaces in $\bC^{n+1}$ which have generic
Levi degeneracies at
$p_0\in M$ and $p_0'\in M'$, respectively. Suppose that the integers
$r$ and conjugacy classes $C_{r,s}$, given by Proposition $1.1.3$, for
$M$ and $M'$ at $p_0$ and $p_0'$, 
respectively, coincide. Then $(M,p_0)$ and
$(M',p_0')$ are biholomorphically equivalent if and only if, for any
choice of $R\in C_{r,s}$ and two
(possibly different) choices of normalization as described in Theorem
$1.1.28$, $(M,p_0)$ and $(M',p_0')$
can be brought to the same normal form.\endproclaim 

\subhead 1.2. The semidefinite case\endsubhead In Proposition 1.1.3, the
partial normal form for a real hypersurface $M$ at a generic
Levi degeneracy $p_0\in M$ is given in terms of a conjugacy class $C_{r,s}$ in
$\Cal S(\bC^{n-1})$. In order to obtain a more explicit partial normal
form, we must distinguish a unique representative in each conjugacy
class. In this paper, we
shall only address this problem in the case where the Levi form
at $p_0$ is semidefinite, i.e. $r=n-1$ and $s=0$, in which case the
group $\hat U(r,s,\bC)$ reduces to the unitary group $U(\bC^{n-1})$ and a lemma due
to E. Cartan can be applied.  The details are worked out in \S5 below. We state
here the corresponding normal forms, which follow from the
results in \S1.1 above and \S5.

Thus, we assume that $M$ has a generic Levi degeneracy
at $p_0\in M$, and that the Levi form at that point is semidefinite
(i.e. the integer $r$ in Proposition 1.1.3 equals $n-1$). An
immediate consequence of Theorem 5.8 (which in fact treats a slightly
more general case; see Theorem 1.2.10 below) is that there are local
holomorphic coordinates $Z=(z,w)$ as in Proposition 1.1.3 such that $M$
is given near $p_0=(0,0)$ by
$$
\multline
\im
w=\sum_{j=1}^{r}|z^j|^2-\sum_{k=r+1}^{n-1} |z^k|^2+
2\re\!\!\left(\bar
z^{\bar n}\left((z')^\tau D_{n-1}(\lambda)z'+(z^n)^2\right)\right)\\+F(z,\bar
z,\re w),\endmultline
\tag1.2.1
$$
where $F$ is as in Proposition 1.1.3, $\lambda=(\lambda_1,\ldots,\lambda_{n-1})$ is a
uniquely determined vector with $\lambda_1\geq\ldots\geq \lambda_{n-1}\geq 0$ such
that either $\lambda_1=1$ or $\lambda_k=0$ for $k=1,\ldots n-1$. We use here the
notation $D_{n-1}(\lambda)$ 
for the diagonal $(n-1)\times(n-1)$-matrix with $\lambda$ on the diagonal, i.e.
$$
D_{n-1}(\lambda)=\pmatrix 
\lambda_1&0&\ldots&0\\
0&\lambda_2&\ldots&0\\
\vdots&\vdots&\ddots&\vdots\\
0&0&\ldots&\lambda_{n-1}
\endpmatrix.\tag1.2.2
$$
An inspection of the proof of Proposition 1.1.9 shows that the most
general transformation of the form \thetag{1.1.6} preserving regular
and the equation \thetag{1.2.1} is of the form
$$
\aligned
\tilde f'(z,w)& =c^{1/2}Uz'+wB+2i(B^* Uz')Uz'+O(3)\\
\tilde f^n(z,w)& =c^{1/3}z^n+O(2)\\
\tilde g(z,w)& =cw+2ic^{1/2}(B^* Uz')w+O(4),\endaligned\tag1.2.3
$$
where $c>0$, $B\in\bC^{n-1}$ (considered as an $(n-1)\times 1$
matrix), and $U\in U(\bC^{n-1})$, if $\lambda=0$, and 
$$
\aligned
\tilde f'(z,w)& =Az'+wB+2i(B^* Az')Az'+O(3)\\
\tilde f^n(z,w)& =z^n+O(2)\\
\tilde g(z,w)& =w+2i(B^* Az')w+O(4),\endaligned\tag1.2.4
$$
where $B\in\bC^{n-1}$ (considered as an $(n-1)\times 1$
matrix), and $A\in U(\bC^{n-1})\cap O_{D_{n-1}(\lambda)}(\bC^{n-1})$, if
$\lambda\neq 0$. (The group $U(\bC^{n-1})\cap O_{D_{n-1}(\lambda)}(\bC^{n-1})$ is
described in more detail in Lemma 5.24.) Using the corresponding
factorization \thetag{1.1.12} and the the description of the space of
normal forms $\Cal N$ given in \S1.1 with $R=D_{n-1}(\lambda)$, we get the
following result. 

\proclaim{Theorem 1.2.5} Let $M$ be a smooth hypersurface in
$\bC^{n+1}$ given
near $0\in M$ by \thetag{1.2.1}, where $\lambda$ is the invariant
$(n-1)$-vector
described
above. Then, given any choice of
normalization (i.e. a choice of $P$ as described above), there is a
unique formal 
transformation \thetag{1.1.6} with this normalization that transforms
the defining equation \thetag{1.2.1} of $M$ at $0$ to
$$
\im
w=\sum_{k=1}^{n-1}|z^k|^2+2\re\!\!\left(\bar
z^{\bar n}\left(\sum_{k=1}^{n-1}\lambda_k(z^k)^2+(z^n)^2\right)\right)+N(z,\bar
z,\re w).
\tag1.2.6
$$
where $N(z,\bar z,s)\in\Cal N$. \endproclaim
 
Due to the explicit description of the normalization of the
transformation to normal form, we can compute a bound on 
the dimension of the stability group $\Aut(M,p_0)$ of a smooth
hypersurface $M\subset\bC^{n+1}$ at a generic semidefinite Levi
degeneracy $p_0\in M$. Recall that $\Aut(M,p_0)$ is the group of
biholomorphic transformations near $p_0$ that fix $p_0$ and map $M$
into itself. It is a real, finite dimensional Lie group in view of
results from [BER3] (see also [Z] and [BER2] for results in the higher
codimensional case). 

\proclaim{Corollary 1.2.7} Let $M\subset\bC^{n+1}$ be a smooth
hypersurface which has a generic semidefinite Levi degeneracy at
$p_0$. Let $\lambda$ be the invariant appearing in \thetag{1.2.1}. Then, the following hold. 
\roster
\item "(a)" If $\lambda=(0,\ldots, 0)$, then
$$
\dim_\bR\Aut(M,p_0)\leq \frac{(n-1)n(n+1)(n+2)}{3}+3n^2-n+1.\tag1.2.8
$$

\item "(b)" If $\lambda=(1,\lambda_2,\ldots, \lambda_{n-1})$ with $1\geq
\lambda_2\geq\ldots\geq \lambda_{n-1}\geq 0$, then we write $(1,u_2,\ldots, u_k, 0)$
for the distinct values of $(1,\lambda_2,\ldots,\lambda_{n-1})$ and denote by
$(m_1$, $m_2,\ldots, m_k, \mu)$ their multiplicities. (Thus, $\mu$ is
the multiplicity of the value $0$.)
Then, 
$$
\multline
\dim_\bR\Aut(M,p_0)\leq
\frac{(n-1)n(n+1)(n+2)}{3}+2n^2+n-1+\\
\sum_{j=1}^k\frac{m_j(m_j-1)}{2}+\mu^2.
\endmultline\tag1.2.9
$$
\endroster
\endproclaim

The bound in Corollary 1.2.7 is sharper than the bound that
follows from the results in [BER2--3]. The latter bound grows like $n^5$
whereas the former grows like $n^4$ as $n\to\infty$. 
The proof of Corollary 1.2.7 consists of counting the
number of parameters in the normalization of the transformation to
normal form and using the explicit representation of $U(\bC^{n-1})\cap
O_{D_{n-1}(\lambda)}(\bC^{n-1})$ provided by Lemma 5.24. The details are
left to the reader. 

Let us conclude this section by mentioning that Theorem 5.8 (in
combination with Theorem 4.15) yields a
partial (third order) normal form in a more general case than the one
considered above. Indeed, as a consequence of Theorem 5.8, we have the
following result, in which the Levi degeneracy is not assumed to be generic.

\proclaim{Theorem 1.2.10} Let $M\subset\bC^{n+1}$ be a real smooth
hypersurface and $p_0\in M$. Suppose that the Levi form of $M$ at
$p_0$ has rank $n-1$ and is semidefinite, i.e. all nonzero eigenvalues
of the Levi form have the same sign. Then, there are local holomorphic
coordinates $Z=(z,w)\in\bC^n\times \bC$ near $p_0$, vanishing at
$p_0$, such that the defining equation of $M$ is of precisely one of
the following forms.
\roster
\item"(i)" For either $\lambda=(1,\lambda_2,\ldots,\lambda_{n-2},0)$ with $1\geq
\lambda_2\geq\ldots\geq \lambda_{n-2}\geq 0$ or $\lambda=(0,\ldots, 0)$,
$$
\im
w=\sum_{k=1}^{n-1}|z^k|^2+2\re\!\!\left(\bar
z^{\bar n}\left(\sum_{k=1}^{n-1}\lambda_k(z^k)^2+2z^{n-1}z^n\right)\right)+F(z,\bar
z,\re w).
\tag1.2.11
$$
\item"(ii)" For either $\lambda=(1,\lambda_2,\ldots,\lambda_{n-2},\lambda_{n-1})$ with $1\geq
\lambda_2\geq\ldots\geq \lambda_{n-1}\geq 0$ or $\lambda=(0,\ldots, 0)$,
$$
\im
w=\sum_{k=1}^{n-1}|z^k|^2+2\re\!\!\left(\bar
z^{\bar n}\sum_{k=1}^{n-1}\lambda_k(z^k)^2\right)+F(z,\bar
z,\re w).
\tag1.2.12
$$
\item"(iii)" For either
$\lambda=(1,\lambda_2,\ldots,\lambda_{n-2},\lambda_{n-1})$ with $1\geq
\lambda_2\geq\ldots\geq \lambda_{n-1}\geq 0$ or $\lambda=(0,\ldots, 0)$,
$$
\im
w=\sum_{k=1}^{n-1}|z^k|^2+2\re\!\!\left(\bar
z^{\bar n}\left(\sum_{k=1}^{n-1}\lambda_k(z^k)^2+(z^n)^2\right)\right)+F(z,\bar
z,\re w).
\tag1.2.13
$$
\endroster
Above, $F(z,\bar z,\re w)$ denotes a smooth, real valued function
which is $O(4)$ in the weighted coordinate system where $z$ has weight
one and $w$ weight two. 

In the case $n=2$, the same result holds with the following modifications:
The only choice for $\lambda$ in {\rm (i)} is $\lambda=0$. In {\rm (ii)} and {\rm (iii)}, both
$\lambda=1$ and $\lambda=0$ are allowed.
\endproclaim

Before proving the results on normal forms presented in these two
sections, we shall introduce a new sequence of invariant
tensors. This will be done in the more general setting of generic
submanifolds of $\bC^N$.

\heading 2. CR invariant tensors\endheading

Let $M\subset\bC^{N}$ be a real generic smooth 
submanifold of codimension $d$. Denote by $T^cM\subset TM$ the complex tangent bundle to
$M$, by $\Cal V=T^{0,1}M\subset\bC T^cM$ the CR bundle of $M$, by
$T^0M\subset T^*M$ the
characteristic bundle, and by $T'M\subset\bC T^*M$ the bundle defined
at each $p\in M$ as the annihilator of $\Cal V_p$. We denote by $n$
the CR dimension of $M$, i.e. $n=N-d$. We have the
following for any $p\in M$
$$
\dim_\bR T^c_pM=2n,\quad \dim_\bC\Cal V_p=n,\quad \dim_\bR T^0_pM=d,\quad
\dim_\bC T'_pM= n+d.\tag2.1
$$
For a vector bundle $E$
over $M$, we denote by $C^\infty(M,E)$ the smooth sections of $E$. The reader is
referred e.g. to [BER4] or [B] for the basics of CR structures. We
shall consider only local properties of $M$ near some point
$p$. Hence, given a point $p\in M$, we may, and we will, identify $M$
with some small open neighborhood of $p$ in $M$.

For a CR vector field $L$ on $M$, i.e. a smooth section of $\Cal V$, we define
an operator $\Cal T_L$ on the smooth 1-forms on $M$ as follows,
$$
\T_L\omega=\frac{1}{2i}L\lrcorner d\omega,\tag2.2
$$
where $\lrcorner$ denotes the usual contraction operator. We should
point out here that we use the notation $\left<\cdot,\cdot\right>$ for
the pairing between $r$-covectors and $r$-vectors
normalized in such a way that if $e_\alpha$ and $e^\beta$, $\alpha,\beta=1,\ldots
m$, are dual bases for an $m$-dimensional vector space $V$ and its
dual $V^*$, respectively, then
$e_{\alpha_1}\wedge\ldots\wedge e_{\alpha_r}$,
$1\leq\alpha_1<\ldots<\alpha_r\leq m$, and
$e^{\beta_1}\wedge\ldots\wedge e^{\beta_r}$, 
$1\leq\beta_1<\ldots<\beta_r\leq m$, are dual bases for
$\Lambda^r(V)$ and $\Lambda^r(V^*)$, respectively (see [St, Chapter
I.4]). This normalization is more 
convenient for our purposes than the one used in e.g. [H], which differs from the
present one by the factor $r!$, and is identical to the one used in [B]. 

We shall
refer to sections of $T'M$ as $(1,0)$-forms and denote by
$\Omega^{1,0}(M)$ the space of smooth $(1,0)$-forms on $M$. It is not
difficult to see that $\T_L\:\Omega^{1,0}(M)\to\Omega^{1,0}(M)$, for if
$\omega\in\Omega^{1,0}(M)$ then, for any CR vector field $K$, we obtain,
by using the well known identity (see [H, Chapter I.2]),
$$
\left<\T_L\omega,K\right>=\left<d\omega,L\wedge K\right>= L\left(\lpb
\omega,K\rpb \right)-K\left(\lpb \omega,L\rpb \right)-\lpb
\omega,[L,K]\rpb =0,\tag2.2
$$
since $\omega$ is a section of $T'M$, which at each point $p\in M$
annihilates $\Cal V_p$, and the CR bundle $\Cal V$ is involutive
(or, as it is also called, formally integrable). We shall use the
notation $\Cal L(M)\subset\Omega^{1,0}(M)$ for for those smooth
$(1,0)$-forms that are sections of $T^0M$. The forms in $\Cal L(M)$
will also be referred to as {\it characteristic forms}.

Let $p\in M$ and let us define a sequence of increasing subspaces
$$
E_{0}(p)\subset E_{1}(p)\subset\ldots\subset E_{k}(p)\subset\ldots\subset
T'_pM.\tag2.3
$$
as follows. Set $E_{0}(p)=\bC\otimes T^0_pM$ and let $E_{j}(p)$, for
$j=1,2,\ldots$, be the linear span of $(1,0)$-covectors of the form 
$$
(\Cal T_{K_{\bar 1}} \ldots  \Cal
T_{K_{\bar j}}\theta)_p,\tag2.4
$$ 
where the $K_{\bar i}$ range over all CR vector fields on $M$ near $p$ and
$\theta$ ranges over the smooth sections of $T^0M$ near $p$. The reason for putting a
bar on the indices of CR vector fields is to be
able to use the notation of tensor algebra in later sections; recall
that the CR vector fields for an embedded CR 
submanifold are really anti-holomorphic vector fields.

We shall see
later that to compute the subspaces $E_{j}(p)$ it suffices to take the
linear span of the covectors \thetag{2.4} where the CR vector fields
$K_{\bar i}$ range over the elements of any basis for the 
CR vector fields near $p$ and the characteristic forms $\theta$ range over a
basis for the smooth sections of $T^0M$ near $p$. We will also show that $M$ is
finitely nondegenerate (see below and also e.g. [BER4]) at $p$ if and
only $E_{k}(p)=T'_pM$ for some $\bar k$. The
reader should also note that these subspaces are the same as, but
differently indexed than, those defined in [E2]. The present
definition is better suited for the purposes of this paper. 

Let us for a given integer $k\geq 0$ denote by $F_k(p)\subset\bar\Cal
V_p$ the subspace of those $\bar N_p\in\bar\Cal V_p$ that annihilate
$E_{k}(0)$, i.e.
$$
F_k(p)=E_{k}(p)^\perp\cap\bar\Cal V_p.\tag2.5
$$
Thus, for $k=0$ we have $F_0(p)=\bar \Cal V_p$. Let $\Cal F_k(M)\subset
C^\infty(M,\bar\Cal V)$ denote the space of those 
sections that take values in $F_k(p)$ at $p$. Note that $\Cal F_k(M)$
is a $C^\infty(M)$-submodule of $C^\infty(M,\bar \Cal V)$. Consider
the following diagram, for integers $j\geq 1$,
$$
\CD
\oversetbrace\text{\rm $j$ times} \to {C^\infty(M,\Cal
V)\times\ldots\times C^\infty(M,\Cal V)}\times \Cal F_{j-1}(M)\times\Cal L(M) @> G_j>> \bC\\
@V e_jVV\\
\undersetbrace\text{\rm $j$ times} \to {\Cal
V_p\times\ldots\times\Cal V_p}\times F_{j-1}(p)\times T^0_pM
\endCD,
\tag2.6
$$
where $e_j$ is the evaluation map at $p$ and $G_j$ is the mapping 
$$
(K_{\bar 1},\ldots,K_{\bar j},\bar N,\theta)\mapsto \lpb\Cal
T_{K_{\bar 1}} \ldots  \Cal
T_{K_{\bar j}}\theta,\bar N\rpb_p,\tag 2.7
$$
We would like to have a multi-linear map 
$$
\psi_j\:\undersetbrace\text{\rm $j$ times} \to {\Cal
V_p\times\ldots\times\Cal V_p}\times F_{j-1}(p)\times T^0_pM\to\bC\tag2.8
$$
that makes the diagram \thetag{2.6} commute. Such a multi-linear map
would, by definition, be an invariant of the CR structure $(M,\Cal V)$ (and hence also
a biholomorphic invariant for the generic submanifold
$M\subset\bC^N$ at $p\in M$).

One of the main results is the following.
\proclaim{Theorem 2.9} For each positive integer $j$, there exists a unique multi-linear mapping 
\thetag{2.8} which makes the diagram
\thetag{2.6} commute. The multi-linear mapping \thetag{2.8}, for each
$j$, is symmetric with respect to permutations of the first $j$
variables.\endproclaim

The multi-linear map $\psi_j$ can also be identified with a tensor 
$$
\psi_j\in \undersetbrace \text{\rm $j$ times}\to {\Cal
V_p^*\otimes\ldots\otimes \Cal
V_p^*}\otimes F_{j-1}(p)^*\otimes (T^0_pM)^*.\tag2.10
$$
Before proving Theorem 2.9, let us make a few remarks.
\remark{Remark $2.11$}
\roster
\item"(i)" For $j=1$ and a fixed characteristic covector $\theta_p\in
T^0_pM$, the Hermitian form $\Cal V_p\times\Cal V_p\to \bC$ defined by
$(L,K)\mapsto \psi_0(L,\bar K,\theta_p)$ coincides with the Levi form
of $M$ at the point $p$ and the characteristic covector $\theta_p$. 

\item"(ii)" As mentioned above and as will be proved below, $M$ is
finitely nondegenerate at $p$ if and only if $E_{k}(p)=T'_pM$ for some
$k$. If $M$ is finitely nondegenerate at $p$, then it is called
$k$-nondegenerate at $p$ if $k$ is the smallest integer for which
$E_{k}(p)=T'_pM$. It follows that for a $k$-nondegenerate CR manifold the tensors
$\psi_j$, $j\geq k+1$, are trivial, since $F_j(p)=\{0\}$. Hence, if
e.g. $M$ is a Levi nondegenerate hypersurface (which is the 
same as a 1-nondegenerate hypersurface), then the only non-trivial
invariant tensor produced by Theorem 2.9 is the Levi form of $M$ at
$p$. 

\item"(iii)" The invariant tensors provide obstructions for two
generic submanifolds $M,M'\subset\bC^N$ of codimension $d$ to be
biholomorphically equivalent at given points $p\in M$, $p'\in M'$. The
submanifolds $(M,p)$ and $(M',p')$ cannot be biholomorhically
equivalent unless $\dim F_j(p)=\dim F'_j(p')$ (with the obvious
notation that corresponding object for $M'$ are denoted with a $'$)
and the tensors 
$\psi_j$ and $\psi'_j$ are equivalent (i.e. there are bases in $\Cal
V_p$, $\Cal V'_{p'}$, $F_j(p)$, $F'_j(p')$, $T^0_pM$, and $T^0_{p'}M'$
such that the representations of $\psi_j$ and $\psi'_j$ are equal) for
each $j=1,2,\ldots$. The reader should note, however, 
that the tensors $\psi_j$ do not provide a {\it complete} set of
invariants in the sense that $(M,p)$ and $(M',p')$ are
biholomorphically equivalent if all tensors are equivalent. This is
illustrated e.g. by Theorem 1.1.28, since the normal form given in that
theorem gives a complete set of invariants (by Corollary 1.1.30) and the
invariants coming from the tensors only enter into the second and
third order terms.
\endroster
\endremark\medskip

\demo{Proof of Theorem $2.9$} We claim that for the
multi-linear mapping $\psi_j$ in \thetag{2.8} such
that the diagram \thetag{2.6} commutes to exist, it is necessary and sufficient
that the following statements hold. 
\roster
\item"(a)" For any ${l}\in\{1,2,\ldots, j\}$, $K'=(K_{\bar 1},\ldots,
K_{\bar l})\in(C^\infty(M,\Cal V))^l$, $A,B\in C^\infty(M,\Cal V)$,
$K''=(K_{\bar l+\bar 1},\ldots, K_{\bar j})\in (C^\infty(M,\Cal V))^{j-l-1}$, $a,b\in
C^\infty(M)$, $\bar N\in 
\Cal F(M)$, and $\theta\in\Cal L(M)$, the following identity holds
$$
\multline
G_j(K',aA+bB,K'',\bar
N,\theta)=\\a(p)G_j(K',A,K'',\bar
N,\theta)+b(p)G_j(K',B,K'',\bar N,\theta).\endmultline\tag2.12
$$

\item"(b)" For any $K=(K_{\bar 1},\ldots, K_{\bar j})\in (C^\infty(M,\Cal V))^j$, $a,b\in
C^\infty(M)$, $\bar A,\bar B\in 
\Cal F(M)$, and $\theta\in\Cal L(M)$, the following identity holds
$$
G_j(K,a\bar A+b\bar B,\theta)=a(p)G_j(K,\bar A,\theta)+b(p)
G_j(K,\bar B,\theta).\tag2.13 
$$

\item"(c)" For any $K=(K_{\bar 1},\ldots, K_{\bar j})\in (C^\infty(M,\Cal V))^j$, $\bar
N\in\Cal F(M)$, $a,b\in
C^\infty(M)$, $\xi,\eta\in\Cal L(M)$, the following identity holds
$$
G_j(K,\bar N,a\xi+b\eta)=a(p)G_j(K,\bar N,\xi)+b(p)G_j(K,\bar N,\eta).\tag2.14
$$
\endroster
Indeed, if the mapping $\psi_j$ exists, then the statements (a), (b),
and (c) follow immediately from the diagram \thetag{2.6}. Conversely, 
if the statements (a), (b), and (c) hold, then the mapping $\psi_j$
can be uniquely constructed as follows. Take $L_{\bar 1},\ldots,
L_{\bar n}$ to be
any basis for the CR vector fields near $p$, $\bar N_1,\ldots, \bar
N_k$ to be generators for $\Cal F_{j-1}(M)$ near $p$ (it is easy
to verify that $\Cal F_k(M)$ is finitely generated as a
$C^\infty(M)$-module near $p$), and $\theta^1,\ldots, \theta^d$ to be a basis
for the characteristic forms near $p$. The restrictions of these
sections to the point $p$ span the corresponding vector space over
$\bC$. We then define $\psi_j(L_{\bar i_1},\ldots, L_{\bar i_j},\bar
N_k,\theta^l)$ to be $G_j(L_{\bar i_1},\ldots, L_{\bar i_j},\bar
N_k,\theta^l)$ and extend $\psi_j$ by linearity. The statements (a),
(b), and (c) guarantee that this definition is independent of the
bases and generators chosen and that the diagram \thetag{2.6}
commutes. We leave the details of this verification to the
reader. These arguments also show that the mapping $\psi_j$ is unique
whenever it exists. 

We begin by proving statement (a). Observe first that the mapping
$G_j$ is clearly multi-linear over $\bC$, so that 
$$
\multline
G_j(K',aA+bB,K'',\bar
N,\theta)=\\G_j(K',aA,K'',\bar
N,\theta)+G_j(K',bB,K'',\bar N,\theta).\endmultline\tag2.15
$$
Hence, it suffices to prove that for any $a$ and $A$ as in statement
(a) we have
$$
G_j(K',aA,K'',\bar
N,\theta)=a(p)G_j(K',A,K'',\bar
N,\theta).\tag2.16
$$
Note that, for any CR vector field $L$, any $b\in C^\infty(M)$, and
any $\omega\in\Omega^{1,0}(M)$, we have
$$
\aligned
\Cal T_L(b\omega) &=L\lrcorner
d(b\omega)=L\lrcorner(db\wedge\omega+bd\omega)\\
& =(L\lrcorner db)\omega-(L\lrcorner\omega)db+b \Cal T_L\omega\\
& =(Lb)\omega+b \Cal T_L\omega,\endaligned\tag2.17
$$
since $L\lrcorner\omega=\lpb \omega,L \rpb=0$. A simple inductive
argument using \thetag{2.17} proves that, for $K_{\bar 1},\ldots,
K_{\bar l}$, $a$,
and $A$ as in the statement (a) and
$\omega$ as above, we have 
$$
\Cal T_{K_{\bar 1}} \ldots \Cal T_{K_{\bar l}}  T_{aA}\omega=a\Cal
T_{K_{\bar 1}} \ldots \Cal T_{K_{\bar l}}  T_{A}\omega+\sum_{i=1}^ma_i
\omega^i,\tag2.18
$$
where the $a_i\in C^\infty(M)$ and the $\omega^i$ are of the form
$$
\omega^i=\Cal T_{S_{\bar 1}} \ldots \Cal T_{S_{\bar k}}\omega\tag2.19
$$
for some $k<l$ and $S_{\bar r}\in\{K_{\bar 1},\ldots, K_{\bar l},A\}$. Hence, for any
$\bar N\in\Cal F_{j-1}(M)$, we have, since $l\leq j$,
$$
\lpb \Cal T_{K_{\bar 1}} \ldots \Cal T_{K_{\bar l}}  T_{aA}\omega,\bar
N\rpb_p=a(p)\lpb \Cal T_{K_{\bar 1}} \ldots \Cal T_{K_{\bar l}}  T_{A}\omega,\bar
N\rpb_p,\tag2.20
$$
which proves \thetag{2.16} if we set $\omega=\Cal
T_{K_{\bar l+\bar 1}} \ldots \Cal T_{K_{\bar j}}\theta$. 

Statement (b) is obvious, since we even have
$$
\lpb\Cal T_{K_{\bar 1}} \ldots  \Cal
T_{K_{\bar j}}\theta,a\bar A+b\bar B\rpb=a\lpb\Cal T_{K_{\bar 1}} \ldots  \Cal
T_{K_{\bar j}}\theta,\bar A\rpb+b\lpb\Cal T_{K_{\bar 1}} \ldots  \Cal
T_{K_{\bar j}}\theta,\bar B\rpb.\tag2.21
$$
Finally, statement (c) follows from an argument similar to the one
used to prove (a). We leave the details to the reader.

To prove the symmetry properties, we first prove the following
identity.
\proclaim{Lemma 2.22} For any $(1,0)$-form $\omega$, CR vector fields
$K$, $L$, and any vector field $X$ on $M$, the following holds
$$
\lpb \Cal T_L \Cal T_K\omega,X\rpb-\lpb \Cal T_k \Cal
T_L\omega,X\rpb = [L,K]\lpb\omega,X\rpb +\lpb
\omega,[X,[L,K]]\rpb.\tag2.23
$$
\endproclaim
\demo{Proof} For any $(1,0)$-form $\omega'$, CR vector field $L'$, and any vector field $X$
on $M$, we obtain, using a well known identity,
$$
\aligned
\lpb \Cal T_{L'}\omega',X\rpb &=\lpb d\omega',L'\wedge X\rpb =
L'\lpb \omega',X\rpb-X\lpb \omega',L'\rpb -\lpb \omega',[L',X]\rpb\\
& =L'\lpb \omega',X\rpb -\lpb \omega',[L',X]\rpb,\endaligned\tag2.24
$$
since $\lpb\omega',L'\rpb =0$. Similarly, we obtain
$$
\lpb \Cal T_L \Cal T_K\omega,X\rpb=LK\lpb\omega,X\rpb -
L\lpb\omega,[K,X]\rpb -K\lpb\omega, [L,X]\rpb - \lpb
\omega,[K,[L,X]]\rpb.\tag2.25
$$
It follows that
$$
\lpb \Cal T_L \Cal T_K\omega,X\rpb-\lpb \Cal T_k \Cal
T_L\omega,X\rpb = [L,K]\lpb\omega,X\rpb +\lpb
\omega,[K,[L,X]]-[L,[K,X]]\rpb.\tag2.26
$$
Now, using the Jacobi identity, we have
$$
[K,[L,X]]-[L,[K,X]]=[K,[L,X]]+[X,[L,K]]+[K,[X,L]]=[X,[L,K]],\tag2.27
$$
which completes the proof.\qed 
\enddemo

In particular, Lemma 2.22 implies that $\Cal T_L$ and $\Cal T_K$,
considered as linear maps on $\Omega^{1,0}(M)$, commute if the CR
vector fields $L$ and $K$ commute. It is well known that there exists
a basis of CR vector fields on $M$ near $p$ that commute. Since this
basis can be used in the construction of $\psi_j$, as described in the
beginning of this proof, it follows that $\psi_j$ is symmetric with
respect to permutations of the $j$ first variables. This completes
the proof of Theorem 2.9.\qed
\enddemo

\heading 3. Finitely nondegenerate CR manifolds\endheading

In this section, we relate the invariant tensors defined in section 2
to the notion of finite nondegeneracy. Let $M\subset\bC^N$ be a
generic real smooth submanifold of codimension $d$, $p_0$ a point in $M$, and let
$\rho(Z,\bar Z)=0$, where $\rho=(\rho_1,\ldots, \rho_d)$,  be a
defining equation for $M$ near $p_0$. Let $L_{\bar 1},\ldots, L_{\bar n}$, $n=N-d$,
be a basis for the CR vector fields of $M$ near $p_0$. 
$M$ is called
{\it finitely nondegenerate} at $p_0$ if there exists a non-negative
integer $k$ such that
$$
\Span\left\{L^{\bar J}\left(\frac{\partial\rho_l}{\partial
Z}\right)(p_0,\bar p_0)\:\forall |J|\leq k,\quad l=1,2,\ldots,
d\right\}=\bC^N,\tag3.1
$$
where we use the notation  $\bar J=(\bar J_1,\ldots, \bar J_k)\in\{1,\ldots, n\}^k$,
$|\bar J|=k$, and $L^{\bar J}=L_{\bar J_1}\ldots L_{\bar J_k}$. If $M$ is finitely
non-degenerate at $p_0$ and $k$ is 
the smallest integer for which \thetag{3.1} holds, then $M$ is called
$k$-nondegenerate at $p_0$. The property of being $k$-nondegenerate
is independent of the choice of defining equations, local coordinates,
and bases for the CR vector fields. Moreover, $M$ is 0-nondegenerate
at $p_0$ if and only if it is totally real at $p_0$, and if $M$ is a
hypersurface, then it is 1-nondegenerate at $p_0$ if and only if it is
Levi-nondegenerate. (See e.g. [BER1] or [BER4] for these statements.)

Finite nondegeneracy was introduced
in [BHR] in connection with a regularity problem for CR mappings of
real hypersurfaces. It was further explored in connection with the
study of holomorphic mappings between generic submanifolds and real
hypersurfaces in [BER1--3]. Finite
nondegeneracy is also related to holomorphic nondegeneracy as
introduced in [S1] (see also [S2]) and essential finiteness as
introduced in [BJT]. The reader is referred to the book [BER4] for
further information and history. 

We prove here the following result. Recall from section 2 the definition of the
subspaces $E_j(p_0)\subset T'_{p_0}M$. 
\proclaim{Theorem 3.2} Let $M\subset\bC^N$ be a generic real
submanifold and $p_0\in M$. Then, $M$ is $k$-nondegenerate at $p_0$ if
and only if $E_k(p_0)=T'_{p_0}M$ and $E_{k-1}(p_0)\subsetneq
T'_{p_0}M$. \endproclaim

Before proving Theorem 3.2, we shall show that the space
$E_k(p_0)$ can be computed in a slightly simpler way than in the
definition given in section 2. Let
$L_{\bar 1},\ldots, L_{\bar n}$ be a basis for the CR vector fields on $M$ near
$p_0$, and $\theta^1,\ldots,\theta^d$ a basis for the characteristic
forms near $p_0$. We shall use the notation $\Cal T^{j}=\Cal
T_{L_{\bar j}}$
and, as above for $J=(J_1,\ldots, J_k)\in\{1,2,\ldots,n\}^k$, we denote by
$$
\Cal T^{J}=\Cal T^{J_1}\circ\ldots\circ\Cal T^{J_k}.\tag3.3
$$

\proclaim{Proposition 3.4} For any nonnegative integer $j$, the
following holds
$$
E_j(p_0)=\Span\{(\Cal T^{J}\theta^l)_{p_0}\:\forall |J|\leq j,\quad
l=1,2,\ldots, d\}.\tag3.5
$$
\endproclaim

\demo{Proof} Observe that the right hand side of \thetag{3.5} is contained 
in $E_j(p_0)$ for any nonnegative $j$. Let $K_{\bar 1},\ldots, K_{\bar
j}$ be
arbitrary CR vector 
fields and $\theta$ an arbitary characteristic form. Since
$L_{\bar 1},\ldots, L_{\bar n}$ and $\theta^1,\ldots, \theta^d$ form bases for the
CR vector fields and the characteristic forms, respectively, near
$p_0$, we have, for $l=1,\ldots, j$,
$$
K_{\bar l}=\sum_{m=1}^na^{\bar m}_{\bar l}L_{\bar m},\quad
\theta=\sum_{i=1}^db_i\theta^i,\tag3.6
$$
for some $a^{\bar m}_{\bar l},b_i\in C^\infty(M)$. The fact that
$(\Cal T_{K_{\bar 1}}\ldots\Cal T_{K_{\bar j}}\theta)_{p_0}$ is contained in the
right hand side of \thetag{3.5} now follows from \thetag{2.17} and
\thetag{2.18}. \qed\enddemo

\demo{Proof of Theorem $3.2$} For a generic submanifold $M\subset\bC^N$
with defining functions $\rho=(\rho_1,\ldots, \rho_d)$ near $p_0\in
M$, we may take $\theta^j=2i\partial \rho_j$, for $j=1,\ldots, d$, as
a basis for the characteristic forms near $p_0$. Observe that each
$\theta^j$ is real on $M$, since $\partial\rho_j+\bar\partial\rho_j=0$
when restricted to $M$. Let $L_{\bar 1},\ldots, L_{\bar n}$ be a basis for the CR
vector fields of $M$ near $p_0$. In the coordinates $Z$ of the ambient
space, we may write
$$
L_{\bar k}=\sum_{l=1}^Na^{\bar l}_{\bar k}(Z,\bar Z)\frac{\partial}{\partial \bar
Z^{\bar l}},\tag3.7
$$
and
$$
\theta^j= 2i\partial\rho_j=2i\sum_{ l=1}^N\frac{\partial\rho_j}{\partial
Z^ l}dZ^ l.\tag 3.8
$$
Hence, using the notation of Proposition 3.4, we have
$$
\Cal T^{k}\theta^j=\sum_{ l=1}^N L_{\bar k}\left(\frac{\partial\rho_j}{\partial
Z^ l}\right)dZ^ l.\tag3.9
$$
Repeating this argument, we obtain
$$
\Cal T^{J}\theta^j=\sum_{l=1}^N L^{\bar J}\left(\frac{\partial\rho_j}{\partial
Z^ l}\right)dZ^ l.\tag3.10
$$
Since we have $(\Cal T^{J}\theta^j)_{p_0}\in T'_{p_0}M$ and since the dimension of
$T'_{p_0}M$ equals $n+d=N$, the conclusion of Theorem 3.2
follows from Proposition 3.4.\qed\enddemo

\heading 4. The third order invariants and a partial normal
form for real
hypersurfaces\endheading 

We shall show that the second and third order tensors $\psi_2$,
$\psi_3$ form a complete set of third order invariants (in a sense
that will be made more precise in Theorem 4.15 below) for real
hypersurfaces. This will be the first step in the proof of Proposition
1.1.3.

Let $M\subset\bC^{n+1}$ be a real
smooth hypersurface. Let 
$L_{\bar 1},\ldots, L_{\bar n}$ be a basis for the CR vector fields on $M$ near some
distinguished point $p\in M$ and $\theta$ a non-zero characteristic
form near $p$. Set
$L_{\alpha}=\overline{L_{\bar\alpha}} $. Denote by $g_{\bar \alpha \beta}$ the components 
of the tensor $\psi_2$ at $p$, which is just the Levi form of $M$ at
that 
point, relative to the bases $L_{\bar 1,p},\ldots, L_{\bar n,p}$ of $\Cal V_p$,
$L_{1,p},\ldots,L_{n,p}$ of $F_0(p)=\bar \Cal V_p$, and
$\theta_p$ of $T^0_pM$, i.e. 
$$
g_{\bar\alpha\beta}=\lpb \Cal T_{L_{\bar\alpha}}\theta, L_{\beta}\rpb_p,\tag4.1
$$
for $\bar\alpha,\beta=1,\ldots, n$. A change of bases
$$
L_{\bar\gamma,p}=b^{\bar\alpha}_{\bar\gamma} L'_{\bar\alpha,p},\quad
\theta_p=a\theta'_p,\tag4.2 
$$
where we use the usual summation convention to raise and lower
indices, yields the transformation 
rule
$$
g'_{\bar\alpha\beta}=a b^{\bar\gamma}_{\bar\alpha}
b^{\nu}_{\beta}g_{\bar\gamma\nu},\tag4.3  
$$
where $b^{\nu}_{\gamma}=\overline{b^{\bar\nu}_{\bar\gamma}}$. By a
suitable choice of bases above, we may assume that the Levi form
of $M$ at $p$ is diagonal with diagonal elements in $\{-1,0,1\}$, i.e.
$$
g_{\bar\alpha\beta}=\epsilon_{\beta}\delta_{\bar\alpha\beta},\tag4.4
$$
where $\delta_{\bar\alpha\beta}$ is the Kronecker symbol and
$$
\epsilon_\beta=\left\{\aligned 1 &,\quad \beta=1,\ldots, r,\\
-1 &,\quad \beta=r+1,\ldots, r+s,\\ 0 &,\quad \beta=r+s+1,\ldots,
n.\endaligned \right.\tag 4.5
$$
We shall assume here that $r+s<n$, so that
$M$ is Levi degenerate at $p$. (The rank of the Levi form at $p$ is $r+s$.)
Now, denote by $h_{\bar\alpha\bar\beta\gamma}$ the
components of the third order tensor $\psi_3$ at $p$, i.e.
$$
h_{\bar\alpha\bar\beta\gamma}=\lpb \Cal T_{L_{\bar\alpha}}\Cal
T_{L_{\bar\beta}}\theta, L_{\gamma}\rpb_p,\tag4.6
$$
where $\bar\alpha,\bar\beta=1,\ldots, n$ and $\gamma=r+s, \ldots,n$.
We then obtain the transformation rule
$$
h'_{\bar\alpha\bar\beta \gamma}=a b^{\bar\sigma}_{\bar\alpha}
b^{\bar\mu}_{\bar\beta} b^{\nu}_{\gamma}h_{\bar\sigma\bar\mu\nu}.\tag4.7
$$

It is well known (and not difficult to see) that we may choose
coordinates $Z=(z,w)=(z^1,\ldots, z^n,w)\in\bC^{n+1}$ near $p\in M$,
vanishing at $p$, such
that $M$ is defined near $p=0$ by the equation $\rho(Z,\bar Z)=0$, where
$$
\rho(Z,\bar Z)=-\im w +
g'_{\bar\alpha\beta}\bar z^{\bar\alpha} z^{\beta}+
2\re\left(k'_{\bar\alpha\bar\beta\nu}\bar z^{\bar\alpha} 
\bar z^{\bar\beta}z^{\nu}\right)+R'(z,\bar z,\re w)\tag4.8
$$
for some $g'_{\bar\alpha\beta}, k'_{\bar\alpha\bar\beta\nu}\in
\bC$ with $\bar\alpha,\beta,\bar\beta,\nu=1,\ldots, n$; here, $R'(z,\bar z,s)$ is a 
real-valued function that vanishes to weighted
order $4$ at $0$ in the 
weighted coordinate system where $z$, $\bar z$  have weight one and $s$
has weight two (or higher if the Levi form at $p$ is 0). For the
embedded hypersurface defined by the function 
\thetag{4.8}, we may take as a basis for the CR vector fields
$$
L'_{\bar \alpha}=\frac{\partial}{\partial
\bar z^{\bar\alpha}}+\lambda_{\bar\alpha}(Z,\bar Z)\frac{\partial}{\partial
\bar w},\quad\bar\alpha=1,\ldots, n,\tag4.9
$$
where $\lambda_{\bar\alpha}(0,0)=0$. We refer the reader e.g. to [BER4, Chapter IV] for
details.  By taking $\theta'=2i\partial\rho$ and using \thetag{3.10},
we find that the tensors $\psi_2$ and $\psi_3$ at $p=0$ relative to
the bases defined by $L'_{\bar\alpha}$, $L'_\beta$, and $\theta'$ are
given by $\psi_2=(g'_{\bar\alpha\beta})$ and
$\psi_3=(h'_{\bar\alpha\bar\beta\gamma})$ with
$$
h'_{\bar\alpha\bar\beta\gamma}=k'_{\bar\alpha\bar\beta\gamma}, \quad
\bar\alpha,\bar\beta=1,\ldots, n,\, \gamma=r+s+1,\ldots, n.\tag4.10
$$
It follows that there is a change of basis \thetag{4.2} such that
\thetag{4.3} and \thetag{4.7} (with $\gamma$ running from $r+s+1$ to
$n$) hold. Such a change of bases corresponds to a linear change of
coordinates of the form
$$
z^\alpha\mapsto b^\beta_\alpha z^\beta,\quad w\mapsto
\frac{1}{a}w\tag4.11
$$
in \thetag{4.8}. Hence, the linear change of coordinates \thetag{4.11}
transforms the defining function in \thetag{4.8} to the form
$$
\multline
\rho(Z,\bar Z)=-\im w +
g_{\bar\alpha\beta}\bar z^{\bar\alpha} z^{\beta}+
2\re\left(k_{\bar\alpha\bar\beta\mu}\bar z^{\bar\alpha} 
\bar z^{\bar\beta}z^{\mu}\right)\\+
2\re\left(h_{\bar\alpha\bar\beta\gamma}\bar z^{\bar\alpha} 
\bar z^{\bar\beta}z^{\gamma}\right)+R(z,\bar z,\re w),\endmultline\tag4.12
$$
where $\bar\alpha,\beta,\bar\beta$ run over $1,\ldots, n$, $\mu$ runs
over $1,\ldots, r+s$, $\gamma$ runs over $r+s+1,\ldots, n$, and
$k_{\bar\alpha\bar\beta\mu}$ are some complex numbers. Next, since
$g_{\bar\alpha\beta}$ is of the form 
\thetag{4.4} with $\epsilon_\beta$ of the form \thetag{4.5}, we
observe that the quadratic change of coordinates
$$
z^{\mu}-\epsilon_{\mu}\overline{k_{\bar\alpha\bar\beta\mu}}
z^\alpha z^\beta\mapsto z^\mu,\tag4.13
$$
for $\mu=1,\ldots, r+s$, yields
the following final form of $\rho(Z,\bar Z)$
$$
\rho(Z,\bar Z)=-\im w +
g_{\bar\alpha\beta}\bar z^{\bar\alpha}z^\beta+2\re\left(
h_{\bar\alpha\bar\beta\gamma}\bar z^{\bar\alpha}
\bar z^{\bar\beta}
z^\gamma\right)+\tilde R(z,\bar z,\re w),\tag4.14
$$
where $\tilde R(z,\bar z,s)$ vanishes of weighted order 4 at 0, the
indices $\bar\alpha,\beta,\bar\beta$ run over $1,\ldots, n$, and the
index $\gamma$ runs over $r+s+1,\ldots, n$. We would like to point out
that a similar form for a real hypersurface was presented by Webster
in [W] (see also Remark 4.17 below).

Hence, we have proved that $\psi_2$ and $\psi_3$ form a complete set of
third order invariants for a real hypersurface $M\subset\bC^{n+1}$ in
the following sense. We use the notation and conventions introduced above.

\proclaim{Theorem 4.15} Let $M\subset\bC^{n+1}$ be a real smooth
hypersurface and $p\in M$. Assume that the signature of
the Levi form of $M$ at $p$ is as described above. Then, there are coordinates
$Z=(z,w)\in\bC^{n+1}$, vanishing at $p$, such that $M$ is defined near
$p=0$ by $\rho(Z,\bar Z)=0$, where $\rho(Z,\bar Z)$ is given by
\thetag{4.14} if and only if there is are bases $L_{\bar 1,p},\ldots,
L_{\bar n,p}$ for $\Cal V_p$, with the corresponding basis $L_{1,p},\ldots,
L_{n,p}$ for $\bar \Cal V_p$, and $\theta_p$ for $T^0_pM$ such that
$$
\psi_2=(g_{\bar\alpha\beta}), \quad
\psi_3=(h_{\bar\alpha\bar\beta\gamma}),\tag4.16
$$
with $\bar\alpha,\beta,\bar\beta=1,\ldots, n$ and
$\gamma=r+s+1,\ldots,n$.\endproclaim 

\remark{Remark $4.17$} In [W], the {\it cubic form} of a real
hypersurface $M\subset\bC^{n+1}$ at a point $p\in M$ was introduced
and shown to be a multi-linear map $\Cal V_p\times\Cal V_p\times
F_1(p)\to\bC$ defined by
$$
q_p(L_p,K_p,\bar N_p)=\lpb \partial\rho,[K,[L,\bar N]]\rpb_p,\tag4.18
$$
where $L$, $K$, and $\bar N$ are vector fields extending
$L_p,K_p\in\Cal V_p$ and $\bar N_p\in F_1(p)$, respectively. A
straightforward calculation, using the formula \thetag{2.24}
repeatedly, shows the 
following relation between the cubic form and the tensor
$\psi_3(\cdot,\cdot,\cdot,\theta)$ (for some
fixed $\theta$ e.g. $\theta=2i\partial\rho$)
$$
\aligned
\psi_3(L_p,K_p,\bar N_p,\theta)-2iq_p(L_p,K_p,\bar N_p)= &\lpb\Cal
T_{L}\Cal T_{K}\theta,\bar N\rpb_p-\lpb \theta,[K,[L,\bar N]]\rpb_p\\=
&\, L\left(\lpb \Cal T_K\theta,\bar 
N\rpb\right)+K\left(\lpb \Cal T_L\theta,\bar N\rpb\right)
\\= &-L\left(\lpb \theta,[K,\bar
N]\rpb\right)-K\left(\lpb \theta,[L,\bar
N]\rpb\right).\endaligned\tag4.19
$$
Nevertheless, the cubic
form and $\psi_3(\cdot,\cdot,\cdot,\theta)$
are in fact equal (possibly modulo some multiplicative constant). This
equivalence 
follows from Theorem 4.15, because it is shown in [W] (using the
notation introduced above) that $M$ can be brought to the form
\thetag{4.14} with 
$$
q_p(L_{\bar\alpha,p},L_{\bar\beta,p},L_{\gamma,p})=\frac{i}{2}
h_{\bar\alpha\bar\beta\gamma},\tag4.20
$$
where $\bar\alpha,\bar\beta,\gamma$ range over the same indices as in
Theorem 4.15. (Thus, Theorem 4.15 is in fact implicit in [W], although
using the cubic form as the third order tensor.)
\endremark\medskip

\heading 5. An explicit computation of the third order tensor in a
special case\endheading

We shall keep the notation and conventions introduced in section 4. 
We would like to compute numerical invariants of the tensor
$\psi_3=(h_{\bar\alpha\bar\beta\gamma})$ under
changes of bases \thetag{4.2} preserving the form \thetag{4.4} of the
second order tensor (the Levi form) $\phi_2=(g_{\bar\alpha\beta})$. 
We shall do this only in the following case, which is a bit more
general than the situation considered in \S1.2. 

We assume that the rank $r+s$ of the Levi form $\psi_2$ at the point equals
$n-1$ and that the Levi form is semidefinite.  (Thus, we do not assume
here that the Levi degeneracy is generic.) We may assume, without 
loss of generality, that the $n-1$ 
nonzero diagonal elements $\epsilon_1,\ldots, \epsilon_{n-1}$ of
$g_{\bar\alpha\beta}$ are +1. We can 
identify the third order tensor $\psi_3$ with a symmetric $n\times n$
matrix $H=(h_{\bar\alpha\bar\beta n})$. 

We associate to each change of
basis in $\Cal V_p$ a matrix $B\in GL(\bC^n)$ by
$B=(b^{\bar\alpha}_{\bar\beta})$. We only consider changes \thetag{4.2} that
preserve the form of $\psi_2$, i.e. such that
$$
aB\tilde IB^*=\tilde I,\tag5.1
$$
where $B^*$ denotes the Hermitian adjoint of $B$ and $\tilde I$ is the
matrix of the Levi form, i.e. in block matrix form 
$$
\tilde I=\pmatrix
I_{n-1}&0\\
0&0
\endpmatrix,\tag5.2
$$
with $I_{n-1}=I_{n-1,0}$ being the $(n-1)\times(n-1)$ identity matrix. It is
easy to see that \thetag{5.1} implies that $B$ must be of the form
$$
B=\pmatrix
V&c\\
0&d
\endpmatrix,\tag5.3
$$
where $c\in\bC^{n-1}$, $d\in \bC$, and $V$ is an
$(n-1)\times(n-1)$-matrix related to $a$ in \thetag{4.2} by
$$
aVV^*=I_{n-1},\tag5.4
$$
i.e. $a>0$ and $\sqrt{a}V$ is a unitary matrix. The transformation
rule \thetag{4.7} for $\psi_3$ becomes 
$$
H'=a\bar dBHB^\tau,\tag5.5
$$
where $B^\tau$ denotes the transpose of $B$. 

Recall that, for a given $(n-1)$-vector 
$$
\lambda=(\lambda_1,\ldots,\lambda_{n-1}),\tag5.6
$$
we denote by $D_{n-1}(\lambda)$ the $(n-1)\times(n-1)$ diagonal matrix
with $\lambda$ on the diagonal (see \thetag{1.2.2}).
We shall also use $e^\tau_{n-1}$ for the $(n-1)$-vector
$$
e_{n-1}^\tau=(0,\ldots,0,1).\tag5.7
$$
If $n=2$, then we take $e^\tau_{n-1}=1$. The main result in this section is
the following, which combined with Theorem 4.15 gives Theorem 1.1.1. We use the matrix 
representations of the second and third order tensors as introduced above.

\proclaim{Theorem 5.8} Let $M\subset\bC^{n+1}$, $n>2$, be a smooth real
hypersurface and $p\in M$. Assume that the Levi form
$g_{\bar\alpha\beta}$ of $M$ at $p$ has rank $n-1$ and is
semidefinite (i.e.\ all nonzero eigenvalues have the same sign). If we
normalize the Levi form  
$g_{\bar\alpha\beta}$ so that its matrix is in the form \thetag{5.2}, then the
matrix $H=(h_{\bar\alpha\bar\beta n})$ of the third order tensor can be
brought to precisely one of the following block matrix forms:
\roster
\item"(i)" For either $\lambda=(1,\lambda_2,\ldots,\lambda_{n-2},0)$ with $1\geq
\lambda_2\geq\ldots\geq \lambda_{n-2}\geq 0$ or $\lambda=(0,\ldots, 0)$,
$$
H'=\pmatrix 
D_{n-1}(\lambda)&e_{n-1}\\
e_{n-1}^\tau&0
\endpmatrix.\tag5.9
$$

\item"(ii)" For either $\lambda=(1,\lambda_2,\ldots,\lambda_{n-2},\lambda_{n-1})$ with $1\geq
\lambda_2\geq\ldots\geq \lambda_{n-1}\geq 0$ or $\lambda=(0,\ldots, 0)$,
$$
H'=\pmatrix 
D_{n-1}(\lambda)&0\\
0&0
\endpmatrix.\tag5.10
$$

\item"(iii)" For either $\lambda=(1,\lambda_2,\ldots,\lambda_{n-2},\lambda_{n-1})$ with $1\geq
\lambda_2\geq\ldots\geq \lambda_{n-1}\geq 0$ or $\lambda=(0,\ldots, 0)$,
$$
H'=\pmatrix 
D_{n-1}(\lambda)&0\\
0&1
\endpmatrix.\tag5.11
$$
\endroster
In the case $n=2$, the same result holds with the following modifications:
The only choice for $\lambda$ in {\rm (i)} is $\lambda=0$. In {\rm (ii)} and {\rm (iii)}, both
$\lambda=1$ and $\lambda=0$ are allowed.
\endproclaim

\remark{Remark $5.12$}\roster
\item"(a)" Note that in the case $n=2$, i.e. in $\bC^3$,
there are only 5 different forms for $H'$ and 
no numerical invariants $\lambda$ (i.e. $\lambda$ is only $1$ or $0$). These 5
forms correspond to the 4 different partial normal forms of type (i)
in [E1, Theorem A] and the case which is not 2-nondegenerate at $p$
(see (b) below). We should point out that the partial normal form of
type (ii) in [E1, Theorem A] corresponds
to an explicit normal form for the third order tensor of a
hypersurface in $\bC^3$ at a point $p$ where the Levi form
vanishes. In this case, there are nontrivial numerical invariant.  
\item"(b)" The only case that corresponds to a hypersurface
$M\subset\bC^{n+1}$ which is not $2$-nondegenerate at $p\in M$ is (ii) with
$\lambda=(0,\ldots, 0)$, i.e. $H'=0$.
\endroster
\endremark\medskip

\demo{Proof of Theorem $5.8$} We assume first that $n>2$. We write the symmetric $n\times
n$-matrix $H$ in block matrix 
form
$$
H=\pmatrix
A&\beta\\
\beta^\tau&\gamma
\endpmatrix,\tag5.13
$$
where $A$ is symmetric $(n-1)\times(n-1)$-matrix, $\beta\in\bC^n$, and
$\gamma\in\bC$. By making a change of bases \thetag{4.2} preserving
the form of $g_{\alpha\bar\beta}$, i.e. the matrix $B$ is of the form
\thetag{5.3} and satisfies \thetag{5.4}, the matrix $H$ transforms
according to the rule \thetag{5.5}. A computation shows that
$$
H'=a\bar d\pmatrix
VAV^\tau +V\beta c^\tau+c\beta^\tau V^\tau+\gamma
cc^\tau&d(V\beta+\gamma c)\\
d(\beta^\tau V^\tau +\gamma c^\tau)& \gamma d^2\endpmatrix.\tag5.14
$$
We shall divide the proof into different cases.

\medskip
\flushpar
{\bf The case $\gamma=0$ and $\beta\neq 0$.} We have
$$
H'=a\bar d\pmatrix
VAV^\tau +V\beta c^\tau+c\beta^\tau V^\tau&dV\beta\\
d\beta^\tau V^\tau&0\endpmatrix.\tag5.15
$$
Let us look for $V$ in the form $V=V_2 V_1$, where $V_1$ is a
unitary matrix such that 
$$
V_1\beta=|\beta| e_{n-1},\tag5.16
$$
with $e_{n-1}$ as defined by \thetag{5.7}. If we write
$A'=V_1AV_1^\tau$, then we have
$$
H'=a\bar d\pmatrix
V_2A'V_2^\tau +|\beta| V_2 e_{n-1} c^\tau+|\beta| c e_{n-1}^\tau V_2^\tau&|\beta|
dV_2 e_{n-1}\\
|\beta| de_{n-1}^\tau V_2^\tau&0\endpmatrix.\tag5.17
$$
If we introduce the vector
$$
p=aV_2^*c\tag5.18
$$
and use the fact that $aV_2V_2^*=I_{n-1}$, then the upper left corner
of $H'$ in \thetag{5.17} can be written
$$
V_2(A'+|\beta|(e_{n-1}p^\tau+p e_{n-1}^\tau))V_2^\tau.\tag5.19
$$
It is easy to check that $p\in\bC^{n-1}$ can be chosen uniquely (which
means that $c$ is determined uniquely as a function of $V_2$ and $a$)
such that $A'+|\beta|(e_{n-1}p^\tau+p e_{n-1}^\tau)$ takes the form 
$$
A'+|\beta|(e_{n-1}p^\tau+p e_{n-1}^\tau)=\pmatrix
E&0\\0&0\endpmatrix,\tag5.20
$$
where $E$ is some symmetric $(n-1)\times(n-1)$-matrix. 
If we write $\tilde V=\sqrt{a} V_2$, then it remains to choose a
unitary matrix $\tilde V$, a positive number $\sqrt{a}$, and a complex
(nonzero) number $d$ so as to normalize the matrix and vector
$$
\bar d\tilde V\pmatrix E&0\\0&0\endpmatrix \tilde V^\tau,\quad
\sqrt{a}|d|^2|\beta|\tilde Ve_{n-1}.\tag5.21
$$
The most general unitary matrix $\tilde V$ satisfying $\tilde
Ve_{n-1}=e_{n-1}$ is of the form
$$
\tilde V=\pmatrix F&0\\0&1\endpmatrix,\tag5.22
$$
where $F$ is a unitary $(n-2)\times(n-2)$-matrix. For such a $\tilde
V$, we get
$$
\tilde V\pmatrix E&0\\0&0\endpmatrix \tilde V^\tau=\pmatrix
FEF^\tau&0\\0&0\endpmatrix.\tag 5.23
$$
At this point we need the following lemma, which is a consequence of
E. Cartan's work on Lie groups
(see [Wa] for a discussion; see also [Si] for the lemma in the present
form). We denote by  
$U(\bC^m)=U^+(m,0,\bC)$ the 
group of unitary transformations in $\bC^m$. We also denote by
$O(\bR^m)$ 
the group of (real) orthogonal transformations in $\bR^m$.

\proclaim{Lemma 5.24} Let $E$ be a symmetric $m\times m$-matrix with
complex matrix elements. Then, there is a unique $m$-vector 
$\lambda=(\lambda_1,\ldots, \lambda_m)$ with $\lambda_1\geq\ldots\geq
\lambda_m\geq 0$ such that 
$$
UEU^\tau=D_m(\lambda),\tag5.25
$$
for some $U\in U(\bC^m)$. In fact, the numbers $\lambda^2_j$ are the eigenvalues of the
positive semidefinite Hermitian matrix $E\bar E$.
Moreover, if $\lambda$ is given as above, and we write $(u_1,\ldots,u_k,0)$ for
the distinct values of $(\lambda_1,\ldots,\lambda_m)$ and $(m_1,\ldots, m_k,\mu)$ for
their multiplicities (e.g. $\mu$ denotes the number of zeros among the
$\lambda_j$), then the subgroup of $U\in U(\bC^m)$ for which
$$
UD_m(\lambda)U^\tau=D_m(\lambda)\tag5.26
$$
consists of all matrices of the form
$$
U=\pmatrix
O_1&0&0&0\\0&\ddots&0&0\\0&0&O_k&0\\0&0&0&V\endpmatrix,\tag5.27
$$
where $O_j\in O(\bR^{m_j})$, $j=1,\ldots, k$, and $V\in
U(\bC^\mu)$. (Observe that $\mu$ could be zero in which case there is
no $V$ in \thetag{5.27}.)
\endproclaim

Now, if we choose $\tilde V$ as in \thetag{5.22} with $F$ chosen such
that $FEF^\tau=D_{n-2}(\lambda')$ for some $(n-2)$-vector $\lambda'$, and then 
$\sqrt{a}>0$ and $d\in\bC\setminus\{0\}$ suitably, 
we obtain $H'$ of the form described by (i) in Theorem 5.8. Also, the
vector $\lambda$, as described in Theorem 5.8 (i), is uniqely determined and it
is clear from the arguments above that 
$H'$ cannot be brought to any of the other forms (ii) or (iii). This
concludes the case $\gamma=0$ and $\beta\neq 0$. 

\medskip
\flushpar
{\bf The case $\gamma=0$ and $\beta= 0$.} In this case, it is clear
from Lemma 5.24 that $H'$ can be brought to the form (ii) (and none of
the forms (i) or (iii)) with $\lambda$, as
described in Theorem 5.8 (ii), uniquely determined.

\medskip
\flushpar
{\bf The case $\gamma\neq 0$.} It is clear from \thetag{5.14} that
we can make the upper right and lower left corner of $H'$ vanish  by
choosing
$$
c=-\frac{1}{\gamma} V\beta.\tag5.28
$$
If we bring the factor $a\bar d$ inside the matrix in \thetag{5.14}
then, with $\tilde V=\sqrt{a}V$ as above, the upper left and lower
right corner of $H'$, respectively, become
$$
\bar d \tilde
V\left(A-\frac{1}{\gamma}\beta\beta^\tau\right)\tilde V^\tau,\quad a\bar d
d^2\gamma.\tag5.29
$$
The equation
$$
a\bar d d^2\gamma=1\tag5.30
$$
determines the argument of $d\in \bC$ uniquely. It also determines the modulus of $d$
uniquely as a function of $a>0$ by
$$
|d|=\frac{1}{|\gamma|^{1/3}a^{1/3}}.\tag5.31
$$
Substituting this into the expression for the upper left corner in
\thetag{5.29} and using Lemma 5.24, we deduce that $H$ can be brought to the form
(iii) (and none of the forms (i) or (ii))  with $\lambda$, as described by
Theorem 5.8 (iii), uniquely determined. This concludes the 
case $\gamma\neq 0$. 

Now, if $n=2$, then a similar, but simpler, argument leads to the
statement concluding Theorem 5.8. \qed 
\enddemo

\heading 6. Proof of Proposition 1.1.3\endheading

We shall keep the notation of \S4--5. Recall that the hypersurface
$M$ is assumed to have a generic Levi degeneracy at $p_0$ at which
point the Levi form has $r$, with $(n-1)/2\leq r\leq n-1$, eigenvalues
of the same sign. Thus, we may assume that the matrix
$(g_{\bar \alpha\beta})$ of the Levi form at $p_0$ equals $I_{r,s}$,
where $I_{r,s}$ is as in \S1.1 and $s=n-1-r$. In view of Theorem 4.15, we must
show that the
matrix $H=(h_{\bar\alpha\bar\beta n})$ of the third order tensor can
be brought, by a change of basis \thetag{4.2} preserving 
the Levi form $(g_{\bar \alpha\beta})=I_{r,s}$, to the form
$$
H=\pmatrix R&0\\0&1\endpmatrix,\tag6.1
$$
for some $R\in S(\bC^{n-1})$, and that, under additional such changes
that also preserve the form \thetag{6.1} of $H$,
the matrix $R$ transforms according to the rule
$$
R'=(cV)^\tau R(cV),\tag 6.2
$$ 
where $c>0$ and $V\in\hat U(r,s,\bC)$ can be chosen arbitrarily. The
most general change \thetag{4.2} preserving the Levi form
$g_{\bar\alpha\beta}$ corresponds to a matrix $B$ as in
\thetag{5.3} with $c\in \bC^{n-1}$, $d\in \bC$, and $\sqrt{|a|}V\in\hat
U(r,s,\bC)$ such that $aVI_{r,s}V^*=I_{r,s}$. If we write $H$ in the
form \thetag{5.13}, then the fact that $M$ has a generic Levi
degeneracy at $p_0$ is expressed by $\gamma\neq 0$, as can be verified
by a straightforward calculation (cf. also [W]). An inspection of the
case $\gamma\neq 0$ in the proof of Theorem 5.8 above shows that $H$
can indeed be brought to the form \thetag{5.32} and $R$ transforms
according to the rule \thetag{6.2}, as desired. This completes the
proof of Proposition 1.1.3.
\qed

\heading 7. Proof of Proposition 1.1.9\endheading

We shall use the notation introduced in \S1.1. Consider a transformation
$$
(z',z^n,w)=(\tilde f'(\tilde z,\tilde w),\tilde f^n(\tilde z,\tilde
w), \tilde g(\tilde z,\tilde w)),\tag7.1
$$
where $(\tilde f',\tilde f^n,\tilde
g)$ is of the form (for convenience, we drop the $\tilde{}$ on the
variables)
$$
\aligned
\tilde f'(z,w) &=Az'+z^nD+wB+z^\tau E z+O(3),\\
\tilde f^n(z,w) &=K^\tau z'+d_n z^n+O(2),\\
\tilde g(z,w) &=cw+2i\left<A'z'+z^nD,\bar B\right>w+O(4),\endaligned\tag7.2
$$
where $A\in GL(\bC^{n-1})$, $D,B,K\in \bC^{n-1}$ (considered as
$(n-1)\times 1$ matrices), $E=(E^\beta)_{1\leq\beta\leq n-1}$ is an
$(n-1)$-vector of $n\times n$ matrices, $d_n\in\bC\setminus\{0\}$,
and $c\in\bR\setminus\{0\}$. This is the most general form of a
transformation that preserves regular coordinates (cf. [E1,
\S5--6]). If we write the formal defining equation of $M$ in the (regular)
coordinates $(\tilde z,\tilde w)$ in complex
form (cf. [BER4] or [E1]), i.e.
$$
\tilde w=\tilde Q(\tilde z,\bar{\tilde z},\bar{\tilde w}),\tag7.3
$$
where $\tilde Q(\tilde z,0,\bar{\tilde w})\equiv \tilde Q(0,\bar{\tilde
z},\bar{\tilde w})\equiv\bar{\tilde w}$, then we obtain, by substituting in
\thetag{1.1.4} and setting $\bar w=0$,
$$
\multline
c\left(1+2i\frac{\left<Az'+z^nD,\bar B\right>}{c}\right)\tilde
Q(z,\bar z,0)=2i\bigg(\left<Az'+z^nD,\bar A\bar z'+\bar z^{\bar n}\bar
D\right>+
\\
\left<B,\bar A\bar z'+\bar z^{\bar n}\bar D\right>\tilde Q(z,\bar z,0)+
\left<z^\tau E z,\bar
A\bar z'+\bar z^{\bar n}\bar D\right>+
\\
(\bar K^\tau\bar z'+\bar d_n\bar z^{\bar n})p_{R}(Az'+z^n D,K^\tau
z'+d_nz^n)\bigg)+\ldots,\endmultline\tag7.4
$$
where $\left<\cdot,\cdot\right>$ and $p_{R}$ are defined by
\thetag{1.1.21} and \thetag{1.1.22}, and where the dots $\ldots$ signify terms that
are either $O(4)$ or of type $(k,l)$ with $l>1$. If the transformation
is to preserve the form \thetag{1.1.4}, then we must have (cf. [E1,\S5])
$$
\tilde Q(z,\bar z,0)=2i(\left<z',\bar z'\right>+2\re\!\!\left (\bar
z^{\bar n} p_{R}(z',z^n)\right)+O(4).\tag7.5
$$
By identifying terms of type $(1,1)$, we deduce that $D=0$ and 
$$
\left<Az',\bar A\bar z'\right>=c\left<z',\bar z'\right>.\tag7.6
$$
Observe that \thetag{7.6} is equivalent to $A^*I_{r,s}A=cI_{r,s}$. 
Identifying terms of type
$(2,1)$ and using 
\thetag{1.1.22}, we also see that $K=0$ and 
$$
\left\{\aligned 
&\frac{\bar d_nd_n^2}{c} =1\\
&\frac{\bar d_n}{c}(A)^\tau RA=R\\
&\left<z^\tau E z,\bar A\bar z'\right>=\frac{2i}{c}\left<Az',\bar
B\right>\left<Az',\bar A\bar z'\right>.\endaligned\right.\tag7.7
$$
The conclusion of Proposition 1.1.19 is now easy to verify. This
completes the proof.\qed

\heading 8. Proof of Theorem 1.1.28\endheading

The proof follows closely the proof
of Theorem B in [E1], which in turn was inspired by the work in
[CM]. The idea is to reduce
the proof to a problem of describing the kernel and range of a certain
linear operator. We shall use the notation introduced in \S1.1.

We write the (formal) defining equation \thetag{1.1.4} of $M$ in the
form
$$
\im w=\left<z',\bar z'\right>+2\re (\bar z^{\bar n}p_R(z))+F(z,\bar
z,\re w),\tag8.1
$$
where $\left<\cdot,\cdot\right>$ is given by \thetag{1.1.21},
$p_R(z)$ by \thetag{1.1.22}, and $F(z,\bar z,s)$ is a formal series
in $\Cal F$ as introduced in \S1.1. 
We subject $M$ to a formal transformation
$$
\tilde z=\tilde f(z,w),\quad \tilde w=\tilde g(z,w),\tag8.2
$$
where $\tilde f=(\tilde f',\tilde f^n)=(\tilde f^1,\ldots,\tilde
f^{n-1},\tilde f^n)$, 
which preserves the form of $M$ modulo terms of weighted
degree at least 4, i.e. the transformed hypersurface $\tilde M$ is
given by a defining equation of the form 
$$
\im \tilde w=\left<\tilde z',\bar {\tilde z'}\right>+2\re (\bar
{\tilde z^{\bar n}}p_R(\tilde z))+\tilde F(\tilde z,\bar
{\tilde z},\re\tilde w),\tag8.3
$$
where $\tilde F(z',\bar z', s')$ is in $\Cal F$. We also require that the new
coordinates are regular for $\tilde M$. Thus, $\tilde f$ and $\tilde
g$ are subjected to the 
restrictions imposed by Proposition 1.1.9. As mentioned in \S1, the most general
transformation of this kind can be factored uniquely as 
$$
(\tilde f(z,w),\tilde g(z,w))=(T\circ P)(z,w),\tag8.4
$$
where $P$ and $T$ are as described in that section. 

To prove Theorem 1.1.28, it suffices to prove that there is a unique
transformation 
$$
T(z,w)=(\hat f(z,w),\hat g(z,w))=(z+f(z,w),w+g(z,w)),\tag 8.5
$$
where $f=(f',f^n)=(f^1,\ldots,f^{n-1},f^n)$,
to normal form (i.e. such that the transformed hypersurface $\tilde M$ is
defined by \thetag{8.3} with $\tilde F\in\Cal N$) such that $f'$ is $O(3)$, $f^n$
is $O(2)$, $g$ is $O(4)$, and such that the constant terms in the
formal series \thetag{1.1.17} vanish. We decompose $(f',f^n,g)$, $F$,
and $\tilde F$ into weighted homogeneous 
parts as follows
$$
\aligned
f'(z,w)=\sum_{\nu=3}^\infty f'_\nu(z,w),\quad f^n(z,w) &=\sum_{\nu=2}^\infty
f^n_\nu(z,w)\quad,\quad g(z,w) =\sum_{\nu=4}^\infty g_\nu(z,w)\\
F(z,\bar z,s) =\sum_{\nu=4}^\infty F_\nu(z,\bar z,s)\quad&,\quad
\tilde F(z,\bar z,s) =\sum_{\nu=4}^\infty \tilde F_\nu(z,\bar z,s).\endaligned
$$
Recall here that $z$ and $\bar z$ are assigned the weight one, $w$ and
$s$ are assigned the weight two, and we say that e.g. $F_\nu(z,\bar
z,s)$ is weighted homogeneous of degree $\nu$ if for all $t>0$
$$
F_\nu(tz,t\bar z,t^2s)=t^\nu F_\nu(z,\bar z,s).
$$
The formal power series $F,\tilde F\in\Cal F$ are related as follows 
$$
\multline
\im \hat g(z,s+i\phi) \equiv\left<\hat f^1(z,s+i\phi),\overline{\hat
f^1(z,s+i\phi)}\right>+
\\
2\re\left(\overline{f^n(z,s+i\phi)}p_R(\hat
f(z,s+i\phi))\right) \,+\tilde F(\hat f(z,s+i\phi),\overline {\hat f}(\bar
z,s-i\phi),\re \hat g(z,s+i\phi)),\endmultline\tag8.6
$$
where
$$
\phi=\phi(z,\bar z,s)=\left<z',\bar z'\right>+2\re (\bar z^{\bar
n}p_R(z,\bar z))+F(z,\bar z,s).\tag8.7
$$
Identifying terms of weighted degree $\nu\geq 4$ we obtain
$$
\multline
F_\nu+\im
g_\nu\equiv \left<z',\overline{f'_{\nu-1}}\right>+\left<f'_{\nu-1},\bar
z'\right>+
\\
(\overline{p_R}+2z^n\bar z^{\bar n})f^n_{\nu-2}+
(p_R+2z^n\bar z^{\bar n})\overline{f^n_{\nu-2}}+\tilde
F_\nu+\ldots,\endmultline\tag8.8 
$$
where
$$
\aligned
F_\nu=F_\nu(z,\bar
z,s)&,\quad \tilde F_\nu=\tilde F_\nu(z,\bar
z,s+i\left<z',\bar z'\right>)\\
\overline{f'_{\nu-1}}=\overline{f'_{\nu-1}}(\bar z,s-i\left<z',\bar z'\right>)&,\quad
f'_{\nu-1}=f'_{\nu-1}(z,s+i\left<z',\bar z'\right>)\quad,\quad\text{\rm etc},
\endaligned\tag8.9
$$
and where  the dots $\ldots$ signify terms that only involve $F_\mu$, $F'_\mu$, $g_\mu$,
$f'_{\mu-1}$, and $f^n_{\mu-2}$ for $\mu<\nu$. We can write this as
$$
\re(ig_\nu+2\left<f'_{\nu-1},\bar
z'\right>+2(\overline{p_R}+2z^n\bar z^{\bar
n})f^n_{\nu-2})=F_\nu-F'_\nu+\ldots.\tag8.10
$$
Let us define the linear operator 
$$
L(f',f^n,g)=\re(ig+2\left<f',\bar
z'\right>+2(\overline{p_R}+2z^n\bar z^{\bar
n})f^n)|_{(z,s+i\left<z',\bar z'\right>)}\tag8.11
$$
from the space $\Cal G$ to the
space $\Cal F$, where $\Cal G$ denotes the space of formal power
series (in $(z,w)$) transformations 
$(f',f^n,g)$ 
such that $f'$ is $O(3)$, $f^n$ is $O(2)$, and $g$ is $O(4)$. Observe
that $L$ maps $(f'_{\nu-1},f^n_{\nu-2}, 
g_\nu)$ to a series that is weighted homogeneous of degree $\nu$. We
note, as in [E1] and [CM], that if we could find subspaces
$$
\Cal G_0\subset\Cal G,\quad\Cal N\subset\Cal F\tag8.12
$$
such that, for any $F\in\Cal F$, the equation
$$
L(f',f^n,g)=F\quad \mod \Cal N\tag8.13
$$
has a unique solution $(f',f^n,g)\in\Cal G_0$ then, given any $F'\in\Cal F$, the
equation \thetag{8.10} would allow us to inductively determine 
the weighted homogeneous parts $F_\nu$ of a normal form $F\in
\Cal N$ and the weighted homogeneous parts
$(f'_{\nu-1},f^n_{\nu-2},g_\nu)$ of the transformation
$(f',f^n,g)\in\Cal G_0$ to normal
form in a unique fashion. (This can
also be formulated as saying that $\Cal G_0$ and $\Cal N$ are
complementary subspaces of the kernel and range of $L$,
respectively). 

Let us therefore define $\Cal G_0\subset\Cal G$
as those $(f',f^n,g)\in\Cal G$ for which the
constant terms in the series \thetag{1.1.17} vanish. 
Thus, the proof of Theorem 1.1.28 will be completed by proving the
following. 

\proclaim{Lemma 8.14} Let $\Cal G_0\subset\Cal G$
be as described above and $\Cal N\subset\Cal F$ as defined in \S$1$. 
Then, for any $F\in\Cal F$, the equation
$$
L(f',f^n,g)=F\quad\mod\Cal N\tag8.15
$$
has a unique solution $(f',f^n,g)\in\Cal G_0$.\endproclaim

\demo{Proof} We shall decompose the equation
\thetag{8.15} according to $(k,l)$-type. We decompose $F\in \Cal F$
as follows
$$
F(z,\bar z,s)=\sum_{k,l}F_{kl}(z,\bar z,s),\tag8.16
$$
where each $F_kl\in\Cal F_{kl}$, i.e. each $F_{kl}$ is in $\Cal F$ and
of type $(k,l)$. We also decompose $(f',f^n,g)\in\Cal G$ as follows
($\beta=1,\ldots, n$)
$$
f^\beta(z,w)=\sum_{k}f^\beta_k(z,w),\quad
g(z,w)=\sum_{k}g_k(z,w),\tag8.17
$$
where $f^\beta_k(z,w)$, $g_k(z,w)$ are homogeneous of degree $k$ in $z$,
e.g.
$$
g_k(tz,w)=t^kg_k(z,w)\,,\quad t>0.\tag8.18
$$
The reader should observe that this redefines e.g. $g_k(z,w)$ which,
previously, denoted the weighted homogeneous part of degree $k$ in
$g(z,w)$. However, in what follows we shall not need the decomposition into weighted
homogeneous terms and, hence, the above notation should cause no
confusion; for the remainder of this section, e.g. $g_k(z,w)$ means the part
of $g(z,w)$ which is homogeneous of degree $k$ in $z$, etc. For brevity, we use
the following notation 
$$
f_w(z,w)=\frac{\partial f}{\partial
w}(z,w),\,\ldots,\,f_{w^m}(z,w)=\frac{\partial^m f}{\partial w^m}(z,w),\,\ldots.\tag8.19
$$
We will use the fact
$$
f(z,s+i\left<z',\bar
z'\right>)=\sum_{m}f_{w^m}(z,s)\frac{(i\left<z',\bar
z'\right>)^m}{m!}.\tag8.20 
$$
We shall identify terms of type $(k,l)$ in \thetag{8.15}. Since the
equation is real, it suffices to consider types where $k\geq l$. Also,
note that for $(k,l)$ such that $\Cal N_{kl}=\Cal F_{kl}$ the
equation \thetag{8.15} is trivially satisfied.

In what follows, we use the notation
$$
F_{kl}=F_{kl}(z,\bar z,s),\quad g_k=g_k(z,s),\quad
\overline{g_k}=\overline{g_k}(\bar z,s) 
\quad,\quad\text{\rm
etc}.
$$
Collecting terms of equal type in
\thetag{8.15}, we obtain the following decoupled systems of differential
equations, for $k\geq 3$,
$$
\left\{
\aligned
\frac{i}{2}g_k &=F_{k0}\\
\left<f'_{k+1},\bar z'\right>+2z^n\bar z^{\bar
n}f^n_k-\frac{\left<z',\bar z'\right>}{2}(g_{k})_{w}
&=F_{k+1,1}\,\mod\Cal N_{k+1,1},
\endaligned
\right.
\tag8.21
$$
and, in addition,
$$
\align
&\left\{
\aligned
p_R\overline{f^n_0}+\frac{i}{2}g_2 &=F_{20}\\
\left<f'_3,\bar z'\right>+2z^n\bar z^{\bar n}f^n_2-i\left<z',\bar
z'\right> p_R(\overline{f^n_{0}})_w-
\frac{\left<z',\bar z'\right>}{2}(g_{2})_{w} &=F_{31}\,\mod \Cal N_{31}\\
i\left<z',\bar z'\right>\left<(f'_3)_w,\bar z'\right>+2i\left<z',\bar
z'\right>z^n\bar z^{\bar n}(f^n_2)_w-\frac{\left<z',\bar
z'\right>^2}{2}p_R(\overline{f^n_0})_{w^2}&+\\\overline{p_R}f^n_4-
\frac{i\left<z',\bar z'\right>^2}{4}(g_2)_{w^2} &=F_{42}
\,\mod\Cal N_{42},\endaligned
\right.
\tag8.22
\\
&\left\{
\aligned
\left<z',\overline{f'_0}\right>+\frac{i}{2}g_1 &=F_{10}\\
-i\left<z',\bar z'\right>\left<z',(\overline{f'_0})_w\right>+\left<f'_2,\bar
z'\right>+2z^n\bar z^{\bar
n}f^n_1&+\\p_R\overline{f^n_1}-\frac{\left<z',\bar
z'\right>}{2}(g_1)_w 
&=F_{21}\,\mod\Cal N_{21}\\
i\left<z',\bar z'\right>\left<(f'_2)_w,\bar
z'\right>-\frac{\left<z',\bar z'\right>^2}{2}
\left<z',(\overline{f'_0})_{w^2}\right>+&\\ 2i\left<z',\bar
z'\right>z^n\bar z^{\bar
n}(f^n_1)_w+\overline{p_R}f^n_3-i\left<z',\bar z'\right>
p_R(\overline{f^n_1})_w&-\\ \frac{i\left<z',\bar z'\right>^2}{4}(g_1)_{w^2} 
&=F_{32}\,\mod\Cal N_{32}\endaligned
\right.
\tag8.23
\\&
\left\{
\aligned
-\im g_0 &=F_{00}\\
2\re(\left<f'_1,\bar z'\right>)+4\re(z^n\bar z^{\bar n}f^n_0)-\left<z',\bar z'\right>\re
(g_0)_w &=F_{11}\,\mod\Cal N_{11}\\
4\left<z',\bar z'\right>\im(z^n\bar z^{\bar
n}(f^n_0)_w)-2\re(\overline{p_R}f^n_2)-&\\2\left<z',\bar
z'\right>\im(\left<(f'_1)_w,\bar z'\right>)+\frac{\left<z',\bar
z'\right>^2}{2}\im (g_0)_{w^2} &=F_{22}\,\mod\Cal 
N_{22}\\
-\left<z',\bar z'\right>^2\re(\left<(f'_1)_{w^2},\bar
z'\right>)-2\left<z',\bar z'\right>^2\re(z^n\bar z^{\bar
n}(f^n_0)_{w^2})-&\\2\left<z',\bar
z'\right>\im(\overline{p_R}(f^n_2)_w)+\frac {\left<z',\bar z'\right>^3}{6}\re
(g_0)_{w^3} &=F_{33}\,\mod\Cal N_{33}.\endaligned
\right.
\tag8.24
\endalign
$$
To show that this system has a 
unique solution
$(f',f^n,g)\in \Cal G_0$, if $\Cal N$ is as defined by \thetag{1.1.25}
and \thetag{1.1.26}, we shall need the following facts. Let $p(z,\bar z)$ be a polynomial
of type $(a,b)$. A direct consequence of a theorem of E. Fischer [F] (see
[S] and [ES]) is the following unique decomposition of any formal series
$F_{kl}\in\Cal F_{kl}$,
$$
F_{kl}=pG_{k-a,l-b}+H_{kl},\tag8.25
$$
where $G_{k-a,l-b}\in\Cal F_{k-a,l-b}$ and $H_{kl}\in\Cal F_{kl}$ with
$$ 
\bar p(\nabla,\bar
\nabla)H_{kl}=0;\tag8.26
$$
here, we use the notation $\bar p(z,\zeta)=\overline{p(\bar z,\bar
\zeta)}$. We shall also need the following lemma, whose proof follows
easily from the decomposition \thetag{8.25} and is left to the
reader. 
\proclaim{Lemma 8.27} Given polynomials $p(z,\bar z)$ and
$q(z,\bar z)$ of type $(a,b)$ and $(c,d)$, respectively, any
$F_{kl}\in\Cal F_{kl}$ can be decomposed in a unique way as follows
$$
F_{kl}=pG^1_{k-a,l-b}+qG^2_{k-c,l-d}+H_{kl},\tag8.28
$$
where $G^1_{k-a,l-b}\in \Cal F_{k-a,l-b}$, $G^2_{k-c,l-d}\in\Cal
F_{k-c,l-d}$, and $H_{kl}\in\Cal F_{kl}$ with
$$
\bar q(\nabla,\bar\nabla)H_{kl}=0,\quad \bar
p(\nabla,\bar\nabla)H_{kl}\in\im S;\tag8.29
$$
here, $S$ is the operator defined by $Su=-\bar p(\nabla,\bar\nabla)(qu)$. 
Moreover, any pair 
$$
(G^1_{k-a,l-b},H_{kl})\in \Cal F_{k-a,l-b}\times \Cal F_{kl}
$$ 
such that \thetag{8.29} holds can occur in such a decomposition \thetag{8.28}.\endproclaim

Now, the system \thetag{8.21--8.24} is very similar to the system
\thetag{9.2.2--9.2.5} in [E1]. To show that there is a
unique solution $(f',f^n,g)\in \Cal G_0$, if $\Cal N$ is as defined by \thetag{1.1.25}
and \thetag{1.1.26}, we proceed more or less exactly as in [E1] and use
the decompositions given by \thetag{8.25} and Lemma 8.27. We leave the
verification to the reader. This completes the proof of Lemma 8.14 and
hence that of Theorem 1.1.28.\qed \enddemo

\enddocument
\end